\documentclass{gtart_h}  

\def\ifplaintex{\expandafter\ifx\csname documentclass\endcsname\relax}

\def\ifplaintex{\expandafter\ifx\csname documentclass\endcsname\relax}

\ifplaintex 
\else
\fi

\expandafter\ifx\csname epsfbox\endcsname\relax\input epsf\fi

\def\gt{{\mathsurround=0pt\it $\cal G\mskip-2mu$eometry \&\ 
$\cal T\!\!$opology}}        %  journal title in recommended style

\def\gtp{{\mathsurround=0pt\it $\cal G\mskip-2mu$eometry \&\ 
$\cal T\!\!$opology $\cal P\!$ublications}}  % GT publications

\def\lognumber#1{\def\thelognumber{#1}}
\def\volumenumber#1{\def\thevolumenumber{#1}}
\def\papernumber#1{\def\thepapernumber{#1}}
\def\volumeyear#1{\def\thevolumeyear{#1}}

\def\pagenumbers#1#2{\def\startpage{#1}\def\finishpage{#2}}
\def\published#1{\def\publishdate{#1}}
\def\proposed#1{\def\theproposer{#1}}
\def\seconded#1{\def\theseconders{#1}}
\def\received#1{\def\receiveddate{#1}}
\def\revised#1{\def\reviseddate{#1}}
\def\accepted#1{\def\accepteddate{#1}}

\def\asciiaddress#1{\def\theasciiaddress{#1}}
\def\asciiemail#1{\def\theasciiemail{#1}}

\long\def\asciiabstract#1{\long\def\theasciiabstract{#1}}

\let\\\par\let\thelognumber\relax
\let\thevolumenumber\relax\let\thepapernumber\relax
\let\thevolumeyear\relax\let\thesamplenumber\relax\let\startpage\relax
\let\finishpage\relax\let\publishdate\relax\let\receiveddate\relax
\let\reviseddate\relax\let\accepteddate\relax\let\theasciititle\relax
\let\theasciiauthors\relax\let\theasciiaddress\relax
\let\theasciiabstract\relax
\let\theasciiemail\relax\let\theshortauthors\relax\let\theshorttitle\relax

\long\def\maketitlep{   % start of definition of \maketitlep

\count0=\startpage

\gt\hfill      %   Journal title (top left) 
\hbox to 77pt{\vbox to 0pt{\vglue -15pt\epsfbox{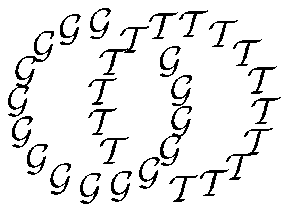}\vss}\hss}
\break
{\small\ifx\thesamplenumber\relax % sample?  
Volume \else Sample
\fi\thevolumenumber\ (\thevolumeyear)
\startpage--\finishpage\nl
Published: \publishdate}
\vglue 0.5truein plus 0.4fil minus 0.1truein

{\parskip=0pt\leftskip 0pt plus 1fil\def\\{\par\smallskip}{\ifplaintex\large
\else\Large\fi\bf\thetitle}\par\medskip}   

\vglue 0pt plus 0.1fil 

{\parskip=0pt\leftskip 0pt plus 1fil\def\\{\par}{\sc\theauthors}
\par\medskip}

\vglue 0pt plus 0.1fil 

{\small\parskip=0pt\let\newline\\
{\leftskip 0pt plus 1fil\def\\{\par}{\sl\theaddress}\par}
\expandafter\ifx\theemail\relax    % email address?
\relax\else\vglue 5pt plus 0.02fil minus 2pt\def\\{\stdspace{\rm 
and}\stdspace} 
\cl{Email:\stdspace\tt\theemail}\fi
\ifx\theurl\relax                  % URL given?
\relax\else\vglue 5pt plus 0.02fil minus 2pt\def\\{\stdspace{\rm 
and}\stdspace}
\cl{URL:\stdspace\tt\theurl}\fi\par}

\vglue 7pt plus 0.3fil minus 3pt

{\bf Abstract}
\vglue 5pt plus 0.1fil minus 2pt

\theabstract

\vglue 7pt plus 0.3fil minus 3pt

{\bf AMS Classification numbers}\quad Primary:\quad \theprimaryclass

Secondary:\quad \thesecondaryclass

\vglue 5pt plus 0.3fil minus 2pt

{\bf Keywords:}\quad \thekeywords

\vglue 10pt plus 0.5fil minus 5pt

{\small  Proposed: \theproposer\hfill Received: \receiveddate\nl
Seconded: \theseconders\hfill 
\ifx\reviseddate\relax                         % paper revised?
Accepted: \accepteddate                        % no
\else
Revised: \reviseddate                          % yes
\fi}
\eject
}       %  end of definition of \maketitlep

\font\phead=cmsl9 scaled 950
\font=cmsl9 scaled 1050
\font\pnum=cmbx10 scaled 913
\font\lnum=cmbx10 
\font\pfoot=cmsl9 scaled 950
\font=cmsl9 scaled 1050
\ifplaintex
\headline{\vbox to 0pt{\vskip -4.5mm\line{\small\phead\ifnum
\count0=\startpage ISSN 1364-0380 (on line)
1465-3060 (printed) \hfill {\pnum\folio}\else\ifodd\count0\def\\{ }% 
\ifx\theshorttitle\relax\thetitle\else\theshorttitle\fi\hfill{\pnum\folio}
\else\def\\{ and }{\pnum\folio}\hfill\ifx\theshortauthors\relax\theauthors
\else\theshortauthors\fi\fi\fi}\vss}}
\footline{\vbox to 0pt{\vglue 0mm\line{\small\pfoot\ifnum\count0=\startpage
\copyright\ \gtp\hfill\else
\gt, Volume \thevolumenumber\ (\thevolumeyear)\hfill\fi}\vss
}}
\else
\makeatletter
\def\@oddhead{{\small\lhead\ifnum\count0=\startpage ISSN 1364-0380 (on line)
1465-3060 (printed) \hfill {\lnum\number\count0}\else\ifodd\count0
\def\\{ }\ifx\theshorttitle\relax \thetitle \else\theshorttitle\fi\hfill
{\lnum\number\count0}\else\def\\{ and }{\lnum\number\count0}
\hfill\ifx\theshortauthors\relax 
\theauthors\else\theshortauthors\fi\fi\fi}}\def\@evenhead{\@oddhead}
\def\@oddfoot{\small\lfoot\ifnum\count0=\startpage\copyright\ \gtp\hfill\else
\gt, Volume \thevolumenumber\ (\thevolumeyear)\hfill\fi}
\def\@evenfoot{\@oddfoot}
\makeatother
\fi

\newwrite\gtoutfile
\long\gdef\makeheadfile{  %%% start of definition of \makeheadfile
{\def\\{, }\def\s{ }
\immediate\openout\gtoutfile head.xxx
\immediate\write\gtoutfile{Proxy-for: \ifx\theasciiauthors\relax
\theauthors\else\theasciiauthors\fi\s<\ifx\theasciiemail\relax\theemail\else\theasciiemail\fi>}
\immediate\write\gtoutfile{\noexpand\\}
\immediate\write\gtoutfile{Authors: \ifx\theasciiauthors\relax
\theauthors\else\theasciiauthors\fi}
{\def\\{ }\immediate\write\gtoutfile{Title: \ifx\theasciititle\relax
\thetitle\else\theasciititle\fi}}
\immediate\write\gtoutfile{Subj-class: GT or SG or MG etc}
\immediate\write\gtoutfile{MSC-class: \theprimaryclass\ifx\thesecondaryclass\relax\else, \thesecondaryclass\fi}
\immediate\write\gtoutfile{Journal-ref: Geom. Topol. \thevolumenumber
(\thevolumeyear) \startpage-\finishpage}
\immediate\write\gtoutfile{Comments: Published by Geometry and Topology at}
\immediate\write\gtoutfile{\s\s http://www.maths.warwick.ac.uk/gt/GTVol\thevolumenumber/paper\thepapernumber.abs.html}
\immediate\write\gtoutfile{\noexpand\\}
\immediate\write\gtoutfile{}
\ifx\theasciiabstract\relax
\immediate\write\gtoutfile{\theabstract}\else
\immediate\write\gtoutfile{\theasciiabstract}\fi
\immediate\write\gtoutfile{}
\immediate\write\gtoutfile{\noexpand\\}
\immediate\write\gtoutfile{}
\immediate\closeout\gtoutfile}}  %%% end of definition of \makeheadfile

\def\maketitlepage{\maketitlep\makeheadfile}
\let\maketitle\maketitlepage

\lognumber{458}
\received{24 May 2004}
\volumenumber{9}\papernumber{18}\volumeyear{2005}
\pagenumbers{757}{811}   
\revised{18 April 2005}
\published{28 April 2005}
\accepted{6 March 2005}
\proposed{Robion Kirby}
\seconded{Ronald Fintushel,  Ronald Stern}

\usepackage{amssymb,amsmath}
\usepackage[curve]{xypic}
\def\S{Section }

\newcommand{\comment}[1]{}
\newbox\mybox
\def\overtag#1#2#3{\setbox\mybox\hbox{$#1$}\hbox to
  0pt{\vbox to 0pt{\vglue-#3\vglue-\ht\mybox\hbox to \wd\mybox
      {\hss$\ss#2$\hss}\vss}\hss}\box\mybox}
\def\undertag#1#2#3{\setbox\mybox\hbox{$#1$}\hbox to 0pt{\vbox to
    0pt{\vglue#3\vglue\ht\mybox\hbox to \wd\mybox
      {\hss$\ss#2$\hss}\vss}\hss}\box\mybox}
\def\lefttag#1#2#3{\hbox to 0pt{\vbox to 0pt{\vss\hbox to
      0pt{\hss$\ss#2$\hskip#3}\vss}}#1}
\def\righttag#1#2#3{\hbox to 0pt{\vbox to 0pt{\vss\hbox to
      0pt{\hskip#3$\ss#2$\hss}\vss}}#1}
\let\ss\scriptstyle

\def\Dot{\lower.2pc\hbox to 2.5pt{\hss$\bullet$\hss}}
\def\Circ{\lower.2pc\hbox to 2.5pt{\hss$\circ$\hss}}
\def\Vdots{\raise5pt\hbox{$\vdots$}}
\def\splicediag#1#2{\xymatrix@R=#1pt@C=#2pt@M=0pt@W=0pt@H=0pt}
\newcommand\lineto{\ar@{-}}
\newcommand\dashto{\ar@{--}}
\newcommand\dotto{\ar@{.}}

\def\omitfootnote#1{}
\newcommand{\interior}[1]{\vphantom{\vbox{\vbox
 {\hbox{$\scriptstyle\circ$}\vskip.3pt}\nointerlineskip\hbox{$#1$}}}
 \vbox{\vbox to 0pt{\vss\hbox{\hskip3pt$\scriptstyle\circ$}
 \vskip.4pt}\nointerlineskip\hbox{$#1$}}}

\newcommand\Brieskorn{Brieskorn}
\newcommand\BCI{BPCI}
\newcommand\BCIs{BPCI's}
\newcommand\CIC{Casson Invariant Conjecture}
\renewcommand\BCI{\Brieskorn{} complete intersection}
\renewcommand\BCIs{\BCI s{}}
\newcommand\sg{\operatorname{sg}}
\newcommand\resgraph{T}

\newcommand\sign{\operatorname{sign}}

\newcommand\Q{{\mathbb Q}}

\newcommand\C{{\mathbb C}}
\newcommand\Z{{\mathbb Z}}
\newcommand\N{{\mathbb N}}
\newcommand\isom{\cong}

\newtheorem{theorem}{Theorem}[section]
\newtheorem{theoremi}{Theorem}
\newtheorem*{theorem*}{Theorem}

\newtheorem*{scholium}{Scholium}
\newtheorem*{fact}{Fact}
\newtheorem{lemma}[theorem]{Lemma}
\newtheorem{proposition}[theorem]{Proposition}
\newtheorem{corollary}[theorem]{Corollary}
\newtheorem{conjecture}{Conjecture}

\newtheorem*{conjecture*}{Conjecture}
\newtheorem*{Cconjecture*}{Casson Invariant Conjecture}
\newtheorem*{stconjecture*}{Splice Type Conjecture}
\theoremstyle{definition}
\newtheorem{definition}[theorem]{Definition}
\newtheorem{example}{Example}

\newtheorem*{examples}{Examples}
\newtheorem{remark}[theorem]{Remark}

\newtheorem*{example*}{Example}
\newtheorem*{remark*}{Remark}
\newtheorem*{notation}{Notation}

\begin{document}
\title{Complex surface singularities with integral\\homology sphere links}
\author{Walter D Neumann\\Jonathan Wahl}
\address{Department of Mathematics, Barnard College, Columbia
  University\\New York, NY 10027, USA}
\address{Department of Mathematics, The University of
North Carolina\\Chapel Hill, NC 27599-3250, USA}
\asciiaddress{Department of Mathematics, Barnard College, Columbia
  University\\New York, NY 10027, USA\\and\\Department of 
Mathematics, The University of
North Carolina\\Chapel Hill, NC 27599-3250, USA}
\gtemail{\mailto{neumann@math.columbia.edu}{\rm\qua 
and\qua}\mailto{jmwahl@email.unc.edu}}
\asciiemail{neumann@math.columbia.edu, jmwahl@email.unc.edu}


\begin{abstract}
While the topological types of {normal} surface singularities with
homology sphere link have been classified, forming a rich class, until
recently little was known about the possible analytic structures.
We proved in \cite{neumann-wahl10} that many of them can be
realized as complete intersection singularities of ``splice type,''
generalizing \Brieskorn{} type.\nl
We show that a normal singularity with homology sphere link is of
splice type if and only if some naturally occurring knots in the
singularity link are themselves links of hypersurface sections of the
singular point.\nl
The Casson Invariant Conjecture (CIC) asserts that for a complete
intersection surface singularity whose link is an integral homology
sphere, the Casson invariant of that link is one-eighth the signature
of the Milnor fiber.  In this paper we prove CIC for a large class of
splice type singularities.\nl
The CIC suggests (and is motivated by the idea) that the Milnor fiber
of a complete intersection singularity with homology sphere link
$\Sigma$ should be a 4--manifold canonically associated to $\Sigma$. We
propose, and verify in a non-trivial case, a stronger conjecture than
the CIC for splice type complete intersections: a precise topological
description of the Milnor fiber.\nl
We also point out recent counterexamples to some overly optimistic
earlier conjectures in \cite{neumann-wahl00} and
\cite{neumann-wahl01}.
\end{abstract}

\asciiabstract{%
While the topological types of {normal} surface singularities with
homology sphere link have been classified, forming a rich class, until
recently little was known about the possible analytic structures.  We
proved in [Geom. Topol. 9(2005) 699-755] that many of them can be
realized as complete intersection singularities of "splice type",
generalizing Brieskorn type. We show that a normal singularity with
homology sphere link is of splice type if and only if some naturally
occurring knots in the singularity link are themselves links of
hypersurface sections of the singular point. The Casson Invariant
Conjecture (CIC) asserts that for a complete intersection surface
singularity whose link is an integral homology sphere, the Casson
invariant of that link is one-eighth the signature of the Milnor
fiber.  In this paper we prove CIC for a large class of splice type
singularities.  The CIC suggests (and is motivated by the idea) that
the Milnor fiber of a complete intersection singularity with homology
sphere link Sigma should be a 4-manifold canonically associated to
Sigma. We propose, and verify in a non-trivial case, a stronger
conjecture than the CIC for splice type complete intersections: a
precise topological description of the Milnor fiber.  We also point
out recent counterexamples to some overly optimistic earlier
conjectures in [Trends in Singularities, Birkhauser (2002) 181--190
and Math. Ann. 326(2003) 75--93].}

{\small\maketitle}

%The topological types of {normal} surface singularities are well
%understood \cite{neumann-transactions}, but it is very rare that
%much is known about the analytic type for given topology. 
In the parallel paper \cite{neumann-wahl10} we give analytic
descriptions in terms of splice diagrams for a wide range of
topologies of singularities, when the link of the singularity is a
$\Q$--homology sphere. The splice diagrams considered there generalize
the original splice diagrams of \cite{eisenbud-neumann, siebenmann} in
that the numerical weights around a node need not be pairwise coprime.
In this paper we restrict to $\Z$--homology sphere links. Our splice
diagrams will thus always have pairwise coprime weights around each
node, and, by \cite{eisenbud-neumann}, the possible links are
classified by their splice diagrams and are obtained by repeatedly
splicing together the links $\Sigma (p_{1},\cdots,p_{n})$ of \BCIs{}
along naturally occurring knots.  Even for this restricted class of
topologies, only in the simplest cases does one know what analytic
properties such singularities might have, eg, being a complete
intersection or Gorenstein.

%explicit analytic descriptions of the singularities were known in only
%the simplest cases.

In \cite{neumann-wahl10}, we describe how ``most''
homology sphere singularity links arise as links of complete
intersection singularities. This occurs when the associated splice
diagram satisfies a certain ``semigroup condition.'' In that case we
give explicit equations, which we call ``splice type,'' generalizing
the \Brieskorn{} complete intersections\omitfootnote{Since others also
  deserved credit we considered the terminology
  ``Brieskorn-\raise1pt\hbox{P}\hglue-1pt\lower1.5pt\hbox{H}am\lower2pt\hbox{m}
  complete intersection.''}. One may think in terms of an operation of
splicing the defining \emph{equations} of two singularities which on
the boundary corresponds to splicing the links.  Specifically, we have
the

\begin{theorem*}{\rm\cite{neumann-wahl10}}\qua Given a homology sphere link
$\Sigma$
whose splice diagram satisfies the semigroup condition, there exists a
complete intersection singularity of splice type whose link is
$\Sigma$.
\end{theorem*}

There is a natural notion of ``higher weight terms'' for a splice type
equation, and, by definition, the result of adding higher weight terms
is still of splice type\footnote{This differs from
  \cite{neumann-wahl00, neumann-wahl01}, where higher order terms were
  not allowed. We now call this ``strict splice type.''} (the effect
on the singularity is always an equisingular deformation). Thus, for
example, the splice type singularities corresponding to one-node
splice diagrams are precisely the \Brieskorn{} complete intersection
singularities with homology sphere link and their higher weight
deformations.

In an earlier paper \cite{neumann-wahl00}, we made the %following
over-optimistic

\begin{stconjecture*}
  Any Gorenstein surface singularity with integral homology sphere
  link is a complete intersection of splice type.
\end{stconjecture*}
\noindent Implicit in this conjecture was a new %and unexpected
necessary
condition (the ``semigroup condition'') on a splice diagram (and hence
on a resolution diagram) in order that it come from a Gorenstein
singularity.  After all, a similar semigroup condition on the value semigroup of
a curve singularity is well known to characterize the Gorenstein ones.
%and a sharpening of this result is in fact a tool in our study.
Further, the conjecture would imply that the \emph{topology} of a
homology sphere link determines a Gorenstein singularity uniquely up
to equisingularity---a kind of ``tautness.''  (Compare
with the equations of plane curve singularities with given
Puiseux pairs.)  Indeed, the
conjecture is true for any singularity $z^n+g(x,y)=0$ with homology sphere
link (Corollary \ref{cor:znfxy}); this is a statement about writing
the irreducible $g(x,y)$ in a certain iterative way.

But a class of examples which may be found in the paper \cite{nemethi et al.} of N\'emethi, Luengo, and Melle-Hernandez, shows this
Conjecture to be false
in this generality:
%, at least with the Gorenstein hypothesis:

\begin{examples}
    {\bf(a)}\qua There exists a Gorenstein singularity, not of splice type, whose
    link is the Brieskorn sphere $\Sigma(2,13,31).$

{\bf (b)}\qua There exists a Gorenstein singularity, not of splice type, whose
    link is a homology sphere but which does not satisfy the semigroup
    conditions.
    \end{examples}

    The above singularities are universal abelian covers of
    ``superisolated'' hypersurface
    singularities.  We do not know in either case how to write down
    equations; in particular, we still know no counterexample to the
    Splice Type Conjecture
    for complete intersections.

We prove our original conjecture under additional assumptions, which clarifies
the situation. A homology
sphere link $\Sigma$ of a normal surface singularity $(X,o)$ has a
number of natural knots, one for each leaf of the splice diagram (or
equivalently, of the resolution graph). For a splice type singularity
these knots are cut out by hyperplane sections.  We prove,
conversely (see Theorem \ref{th:ends} for a more precise version):

\begin{theoremi}\label{thi:ends}
  For a normal surface singularity $(X,o)$ with homology sphere
  link, if all the knots associated to leaves of the splice
  diagram are links of hypersurface sections of\/ $X$, then the
  semigroup condition is fulfilled, and $X$ is a complete intersection
  of splice type.
\end{theoremi}

Our study of singularities with homology sphere link originated in our
conjecture, formulated in \cite{neumann-wahl90}:

\begin{Cconjecture*}
  Let $(X,o)$ be an isolated complete intersection surface singularity
  whose link $\Sigma$ is an integral homology 3--sphere.  Then the Casson
  invariant $\lambda(\Sigma)$ is one-eighth the signature of the
  Milnor fiber of $X$.
\end{Cconjecture*}
\noindent
At the time, we verified the \CIC{} for \Brieskorn{} complete
intersections by direct computation.  It was a challenge to find other
examples, but having done so, the conjecture was verified in these cases,
with the serious work being calculation of the signature.  With the
singularities of splice type we now have an abundance of examples, but
even for these the signature calculations are difficult, and we cannot
verify the \CIC. Still, the
following theorem includes all previously proved cases of the Casson
Invariant Conjecture, except for some cases described by Collin
and Saveliev in \cite{collin-saveliev} (see Remark
\ref{collin-saveliev}).

\begin{theoremi} \label{thi:2node}
  The Casson Invariant Conjecture is true for complete intersection
  singularities of splice type for which the nodes of the splice
  diagram are in a line.
\end{theoremi}

This is proved by reformulating (as in \cite{neumann-wahl90})
the \CIC{} in terms of geometric genus, which is
easier than the signature to compute from defining equations.
\begin{Cconjecture*}[Version 2]   Let $(X,o)$ be a complete
intersection
  surface singularity with integral homology 3--sphere link $\Sigma$.
  Then the Casson invariant $\lambda(\Sigma)$ equals
  $-p_g(X,o)-\frac18C(\Sigma)$, where $C(\Sigma)$ is the characteristic
  number $c_1^2+c_2-1$ of any good resolution of $X$ (this is a
  topological invariant).
\end{Cconjecture*}
This version is equivalent to the previous version by formulas of
Laufer and Durfee (see proof of Theorem \ref{th:equiv}).  This
formulation makes sense for Gorenstein singularities, but is false in
that generality, as seen using some of the examples above.

Assuming the Splice Type Conjecture for complete intersections
(a shaky assumption), one might expect to verify the Casson
Invariant Conjecture by direct calculation with the equations.
But we expect things to go in the
opposite direction: a proof of the Casson Invariant Conjecture
(perhaps symplectic or gauge-theoretic) might help deduce the
form of defining equations.  This happens for instance in the one-node
case: we proved in \cite{neumann-wahl90} that a Gorenstein singularity
$(X,o)$ with link $\Sigma (p_{1},\cdots,p_{n})$ is of splice type,
ie, an equisingular deformation of the corresponding \BCI{}, if
and only if the Casson Invariant Conjecture holds for $X$
(equivalently, $X$ has the same geometric genus as the \BCI).  We
remark that A.  N\'emethi \cite{nemethi} has proved this value of
geometric genus for weakly elliptic singularities, eg, when the link
is $\Sigma(2,3,6k+5)$.

Part of the interest of the Casson Invariant Conjecture is its
suggestion that the Milnor fiber is a ``natural'' 4--manifold which is
attached to its boundary $\Sigma$, and for which the signature
computes the Casson invariant exactly (and not just mod 2).
Specifically, it implies that for a complete intersection singularity
whose link is a homology sphere, analytic invariants like the Milnor
number and geometric genus are determined by the link.  (Such results
are known to be false for general hypersurface singularities.)  Given
the equations of a singularity,
it is relatively easy to calculate the Casson invariant of
the link, but it is extremely hard to calculate the signature of the
Milnor fiber (let alone understand its
topology).

We conjecture a topological construction that, when splicing two
singularities, creates the new Milnor fiber out of the old ones,
extending the operation of splicing on the boundaries (see Conjecture
\ref{conj:milnor fiber}).  This conjecture easily implies
%(and hence motivates)
the Casson Invariant Conjecture for splice type singularities
(Corollary \ref{cor:milnor fiber}).  We succeed in proving it in
a non-trivial case:
\begin{theoremi}\label{thi:fiber}
  For a singularity $z^n+g(x,y)=0$ with homology sphere link, the
  Milnor fiber is formed by the conjectured topological construction.
\end{theoremi}
\noindent Though the Casson Invariant Conjecture for this case
follows, it had already been proven in \cite{neumann-wahl90} (by a
much less conceptual proof), and more recently by Collin and Saveliev
\cite{collin-saveliev} using equivariant Casson invariants and by
N\'emethi and Nicolaescu \cite{nemethi-nicolaescu3} in a more general
context. It is also a special case of Theorem \ref{thi:2node}.

In \cite{neumann-wahl00, neumann-wahl01} we proposed a more
general version of the Splice Type Conjecture: \emph{Any
  $\Q$--Gorenstein surface singularity with $\Q$--homology sphere link
  has as universal abelian cover a complete intersection singularity
  of splice type} (using a more general notion of splice diagram).
Although true surprisingly often, the examples of \cite{nemethi et al.}
mentioned above show this to be false in general, even for
hypersurface singularities.

The converse direction,
that equations of splice type lead to abelian covers of
$\Q$--Gorenstein singularities with expected topological type, is the
main content of \cite{neumann-wahl10} (in particular, as already
mentioned, the equations of splice type of the current paper give
singularities with the expected homology sphere
links). %, so we do not give a proof here.
This paper is nevertheless somewhat transverse to
\cite{neumann-wahl10}, since %in \cite{neumann-wahl01, neumann-wahl10}
we offer there no guess as to the topology or the signature of the
Milnor fiber of the universal abelian cover.  Though
\cite{neumann-wahl90} wondered about a generalization of the Casson
Invariant Conjecture for $\Q$--homology sphere links involving the
Casson--Walker invariant, computations for Seifert fibered rational
homology spheres by Lescop \cite{lescop90, lescop96} showed the naive
generalization fails (see also \cite{collin}).  Lim's result
\cite{lim} suggested looking at a Seiberg--Witten invariant, and a
recent generalization along these lines of the Casson Invariant
Conjecture to $\Q$--Gorenstein $\Q$--homology spheres has been offered
by N\'emethi and Nicolaescu
\cite{nemethi-nicolaescu,nemethi-nicolaescu2,nemethi-nicolaescu3}, but
is now also known to be false in the generality stated (see
\cite{nemethi et al.}).

We offer now a road map to help readers go through this
paper.

Sections \ref{sec:sd} and \ref{sec:equations} are introductory. In
Section \ref{sec:sd}, we review from \cite{eisenbud-neumann} the
definition of splice diagrams and the topological description of
homology sphere links; further details are found in the Appendix
(Section \ref{sec:splicing}), where we also give an improved
description of the relationship between splice diagrams and plumbing
(or resolution) graphs. We also introduce the important ``semigroup
condition.''  In Section \ref{sec:equations} we associate ``splice
type equations'' to any splice diagram with semigroup
condition; this provides a wealth of examples of
complete intersections with homology sphere links.  Modifying the
construction provides familiar equations for complete intersection
monomial curves.

Section \ref{sec:semigroup} develops some theory of
semigroups and monomial curves that is needed in the next two
sections to prove Theorems \ref{thi:ends} and \ref{thi:2node}.
In particular, it includes a new characterization of
complete intersection monomial curves in terms of one-dimensional
analogues of splice type singularities (Theorem \ref{th:cisg} and
its scholium).

Section \ref{sec:semigroup condition} examines the key
property of a splice type singularity: the natural knots in the link
associated to leaves in the splice diagram are obtained by setting a
coordinate equal to 0.  We prove (Theorem \ref{th:ends}, a more
precise version of Theorem \ref{thi:ends}) that
conversely any normal surface singularity with homology sphere link,
and for which the natural knots are hypersurface sections, is in fact
a splice type singularity.  Major use is made of Theorem
\ref{th:cisg} concerning the $\delta$--invariant of certain monomial
curves.  

Section \ref{sec:genus} has as its goal the inductive calculation of
the geometric genus $p_{g}$ for a splice type singularity.  Every node
$v$ of the splice diagram gives a valuation (or weight function) $\nu$
of the singularity; a key result (Theorem \ref{th:assocgraded}) states
that the associated graded ring associated to $\nu$ is an integral
domain, whose normalization is a \Brieskorn{} complete intersection.
Now, $p_{g}$ is the colength of the ``canonical ideal,'' given by
functions for which every $\nu$--weight is at least some explicit
value.  When all the nodes of the splice diagram are on a line, there
is a simultaneous monomial basis for every associated graded (Lemma
\ref{lem:7.5}).  This reduces the calculation of $p_{g}$ in that case
to counting integral lattice points in some region; an induction now
works, yielding the main result, Theorem \ref{th:2n}. Theorem
\ref{thi:2node} is a corollary of this and Theorem \ref{th:equiv} of
the next section. 

The remaining sections \ref{sec:milnor fiber} to \ref{sec:plane
  curves} discuss the Milnor Fiber Conjecture and are largely
independent of the preceding sections.  Section \ref{sec:milnor fiber}
introduces this conjecture, which describes the conjectured topology
of the Milnor fiber of splice type singularities, and which would
imply the \CIC{}.  The discussion leads to Theorem \ref{th:equiv},
which clarifies how the Casson Invariant
Conjecture relates to splicing. This involves the relationship between
signature and geometric genus, and the key is to understand the
behavior of the topological invariant $C(\Delta)$ of the link under
splicing. This is done in Theorem \ref{th:Cdelta}, whose proof, using
numerics of splice diagrams, takes up the following section (Section
\ref{sec:canonical}).

Section \ref{sec:plane curves} verifies the Milnor Fiber Conjecture
for equations of the form $z^{n}=f(x,y)$, by careful topological
construction of the Milnor fiber. This uses a description of plane
curve singularities in terms of splice diagram equations.

{\bf Acknowledgements}\qua 
The conjectures and some results of this paper arose from a visit by
the first author to Duke University, and we thank the Duke Mathematics
Department for its hospitality. We also thank the Max-Planck-Institut
f\"ur Mathematik in Bonn for its hospitality while some of the work on
this paper was done.

The first author's research is supported under NSF grant DMS-0083097
and the second author's under NSA grant MDA904-02-1-0068.

\section{Splice diagrams for integral homology sphere links}\label{sec:sd}
%\section{Equations associated to splice diagrams}

For more details on splicing see the Appendix (Section \ref{sec:splicing}).

Recall that a \emph{splice diagram} is a finite tree with vertices
only of valency 1 (``\emph{leaves}'') or $\ge3$ (``\emph{nodes}'') and
with a collection of integer weights at each node, associated to the
edges departing the node. The following is an example.
$$\splicediag{12}{30}{
\Circ&&&\Circ\\
&\Circ\lineto[ul]_(.25){2}\lineto[dl]^(.25)3
&\Circ\lineto[dr]_(.25){5}\lineto[ur]^(.25){2}
\lineto[l]_(.2){11}_(.8){7}\\
\Circ&&&\Circ
}$$
For an edge connecting two nodes in a splice diagram the \emph{edge
  determinant} is the product of the two weights on the edge minus the
product of the weights adjacent to the edge. Thus, in the above
example, the one edge connecting two nodes has edge determinant
$77-60=17$.

The splice diagrams that classify homology sphere singularity
links satisfy the following conditions on their weights:
\begin{itemize}
\item the weights around a node are positive and pairwise coprime;
\item the weight on an edge ending in a leaf is $>1$;
\item all edge determinants are positive.
\end{itemize}

More general splice diagrams appear for other situations (see,
eg, \cite{eisenbud-neumann, neumann-wahl01, neumann-wahl10}), but we
will only consider splice diagrams satisfying the above conditions
here.

\begin{theorem}{\rm\cite{eisenbud-neumann}}\qua
The homology spheres that are singularity links are in one-one
correspondence with splice diagrams satisfying the above conditions.
\end{theorem}

The splice diagram and resolution diagram for the singularity
determine each other uniquely, and describe how to construct the link
by splicing or by plumbing.  One method to compute the resolution diagram
from the splice diagram is given in \cite{eisenbud-neumann}. We
describe an easier method in the appendix to this paper (Section
\ref{sec:splicing}), where we also recall the topological meaning of
splicing and how to compute the splice diagram from the
resolution diagram for a singularity.

The following notations will be used extensively in this paper.
\begin{notation}
  For a node $v$ and an edge $e$ at $v$, let $d_{ve}$ be the weight on
  $e$ at $v$, and $d_v$ the product of the $d_{ve}$ over all such $e$.
  Let $\Delta_{ve}$ be the subgraph of $\Delta$ cut off from $v$ by
  $e$.  For any pair of vertices $v$ and $w$, let $\ell_{vw}$ (the
  \emph{linking number}) be the product of all the weights adjacent
  to, but not on, the shortest path from $v$ to $w$ in $\Delta$. We
  also consider $\ell'_{vw}$, the same product but excluding weights
  around $v$ and $w$. Thus if $v$ is a node and $w$ is a leaf in
  $\Delta_{ve}$, then
$$\ell_{vw}d_{ve}=\ell'_{vw}d_v.$$
\end{notation}
\begin{definition}[Semigroup Condition]
Let $\Delta$ be a splice diagram. We say $\Delta$
satisfies the
  \emph{semigroup condition} if, for each node $v$ and adjacent edge
  $e$, the edge-weight $d_{ve}$ is in the semigroup
$$\N\langle\ell'_{vw}:w \text{ a leaf of $\Delta$ in }
  \Delta_{ve}\rangle\,.$$
Equivalently, the product $d_v$ of the edge-weights adjacent to $v$ is
in the semigroup
$$\N\langle\ell_{vw}:w \text{ a leaf of $\Delta$ in }
  \Delta_{ve}\rangle\,.$$
\end{definition}
For instance, in the two-node splice diagram above, let $v$ be the
leftmost node and $w$ the upper right hand leaf.  Then $\ell_{vw}$
equals $2\cdot3\cdot5$, while $\ell'_{vw}=5$; the semigroup condition
is satisfied at that node since $7$ is in the semigroup generated by
$2$ and $5$.

If a splice diagram satisfies the semigroup condition, we will write
down complete intersection equations that give a singularity with the
given link. We know of no counterexample to the following optimistic
conjecture mentioned in \cite{neumann-wahl01}.

%The following is a special case of
%Conjecture 2 of \cite{neumann-wahl01}.
%\begin{conjecture}[Gorenstein implies Semigroup Condition]\label{conj:sg}
%If a surface singularity with homology sphere link is Gorenstein, then
%its splice diagram satisfies the semigroup condition.
%\end{conjecture}

\begin{conjecture}[Complete Intersection implies
Semigroup Condition]\label{conj:sg}
  If a surface singularity with homology sphere link is a complete
  intersection, then its splice diagram satisfies the semigroup
  condition.
\end{conjecture}

For example, consider the splice diagram
$$
\splicediag{6}{30}{\\
&\Circ&&&\Circ\\
\Delta\quad=&&\Circ\lineto[ul]_(.25){p}\lineto[dl]^(.25){q}
&\Circ\lineto[dr]_(.25){q'}\lineto[ur]^(.25){p'}
\lineto[l]_(.2){r'}_(.8){r}\\
&\Circ&&&\Circ
}$$
with $p,q,r$ and
$p',q',r'$ pairwise coprime triples of positive integers satisfying
$rr'>pqp'q'$.  Then $\Delta$ satisfies the semigroup condition if and only if
$$r\in \N\langle p',q'\rangle\quad\text{and}\quad r'\in \N\langle
p,q\rangle.$$
(Note $r$ is automatically in the semigroup $\N\langle p',q'\rangle$
if it is greater than or equal to the conductor $(p'-1)(q'-1)$.)
In particular, the resolution diagram
$$
\xymatrix@R=6pt@C=24pt@M=0pt@W=0pt@H=0pt{
\\
&\overtag{\Circ}{-2}{8pt}&&&\overtag{\Circ}{-2}{8pt}\\
{\resgraph\quad=}&&\overtag{\Circ}{-7}{8pt}\lineto[ul]\lineto[dl]\lineto[r]&
\overtag{\Circ}{-1}{8pt}\lineto[ur]\lineto[dr]\lineto[l]&\\
&\overtag{\Circ}{-3}{8pt}&&&\overtag{\Circ}{-3}{8pt}}
$$
gives the splice diagram
$$
\splicediag{6}{30}{\\
&\Circ&&&\Circ\\
\Delta\quad=&&\Circ\lineto[ul]_(.25){2}\lineto[dl]^(.25)3
&\Circ\lineto[dr]_(.25){3}\lineto[ur]^(.25)2
\lineto[l]_(.2){37}_(.8){1}\\
&\Circ&&&\Circ
}$$
which does not satisfy the semigroup condition, since $1$ is not in
the semigroup generated by $2$ and $3$. We would therefore expect that
there is no
%Gorenstein
complete intersection singularity with this resolution.

\section{Equations associated to a splice diagram}\label{sec:equations}

Let $\Delta$ be a splice diagram satisfying the semigroup condition.
We will write down a system of complete intersection equations that
give a singularity with the corresponding link.  Associate a variable
$z_w$ to each leaf $w$ of the splice diagram. To each node $v$ of the
splice diagram, we will associate $(\delta_v-2)$ equations, where
$\delta_v$ is the valency of the node.  If $n$ is the number of
leaves, then it is easy to check that $n-2 = \sum(\delta_v-2)$ (summed
over the nodes of $\Delta$), so this will give the right number of
equations.

Fix a node $v$. For each leaf $w$ we give the variable $z_w$ weight
$\ell_{vw}$ (we call this the \emph{$v$--weight} of $z_w$). For each
edge $e$ at $v$ the semigroup condition lets us write
\begin{equation}\label{eq:sg}
    d_{v}=\sum_w\alpha_{vw}\ell_{vw}\,,\quad\text{sum over the leaves
       $w$ of $\Delta$ in
       $\Delta_{ve}$, with $\alpha_{vw}\in\N$. }
\end{equation}
Equivalently,
\begin{equation}\label{eq:sg'}
d_{ve}=\sum_w\alpha_{vw}\ell'_{vw}\,,\quad\text{sum over the leaves
       $w$ of $\Delta$ in $\Delta_{ve}$}
\end{equation}
%This allows us to
We define an \emph{admissible monomial} (associated to
the edge $e$ at the node $v$) to be a monomial $\prod_w
z_w^{\alpha_{vw}}$, the product over leaves $w$ in $\Delta_{ve}$, with
exponents satisfying the above equations.  Thus an admissible monomial
$M_{ve}$ associated to $v$ has total $v$--weight $d_{v}$ (and depends
on the choice of $\alpha_{vw}$).

Next, choose one admissible monomial $M_{ve}$ for each edge at $v$ and
consider $\delta_v-2$ equations associated to $v$ by equating to $0$
some $\C$--linear combinations of these monomials:
$$\sum_e a_{ie}M_{ve}=0,\quad i=1,\dots,\delta_v-2.$$
Repeating for
all nodes, we get a total of $n-2$ equations. If the coefficients
$a_{ie}$ of the equations are ``sufficiently general,'' we say that
the resulting system of $n-2$ equations is of \emph{strict splice type}.

\emph{Sufficiently general} simply means that for every $v$, all
maximal minors of the $(\delta_v-2)\times\delta_v$ matrix $(a_{ie})$
of coefficients should be non-singular.  By applying row operations to
such a matrix (taking linear combinations of the equations) one can
always put the $(\delta_v-2)\times\delta_v$ coefficient matrix in the
form
$$
\begin{pmatrix}
  1&0&\dots&0&a_1&b_1\\
0&1&\dots&0&a_2&b_2\\
\vdots&\vdots&&\vdots&\vdots&\vdots\\
0&0&\dots&1&a_{\delta_v-2}&b_{\delta_v-2}
\end{pmatrix}
$$
so we will often assume we have done so. In this way, the defining
equations are sums of three monomials.  The ``sufficiently general''
condition is then $a_ib_j-a_jb_i\ne0$ for all $i\ne j$, and all $a_i$
and $b_i$ nonzero.

\begin{example}
  Assume $\Delta$ has one node, of valency $n$.  There is no semigroup
  condition. There is only one admissible monomial for each edge,
  namely $z_j^{d_j}$, where $d_j$ is the weight on the edge.  Our
  equations are thus of \Brieskorn{} type:
  $$\sum_{j=1}^na_{ij}z_j^{d_j}=0,\quad i=0,\dots,n-2\,.$$
  The
  ``sufficiently general'' condition is then the well-known condition
  (due to H. Hamm \cite{hamm}) for the system of $n-2$ equations to
  have an isolated singularity.  Thus, for a splice diagram with one
  node, ``strict splice type'' is equivalent to isolated \BCI.
\end{example}
\begin{example}
  For the $\Delta$ of the example at the start of Section \ref{sec:sd}
we associate variables
  $z_1,\dots,z_4$ to the leaves as follows:
  $$
  \splicediag{6}{10}{\\
    &&z_1&\Circ&&&&&&&&&\Circ&z_4\\
    \Delta\quad=&&&&&&\Circ\lineto[ulll]_(.25){2}\lineto[dlll]^(.25)3
    &&&\Circ\lineto[drrr]_(.25){5}\lineto[urrr]^(.25)2
    \lineto[lll]_(.2){11}_(.8){7}\\
    &&z_2&\Circ&&&&&&&&&\Circ&z_3 }\qquad$$
  The admissible monomials
  for the left node are $z_1^2$, $z_2^3$, and $z_3z_4$. The admissible
  monomials for the right node are $z_3^5$, $z_4^2$, and $z_1z_2^4$
  or $z_1^3z_2$ (since $11=\alpha \cdot 3+\beta \cdot 2$ has
  solutions $(1,4)$ and $(3,1)$).  Thus the system of equations might
  be
$$
%\begin{equation}\label{eq1}
\begin{array}{r}
  z_1^2+z_2^3+z_3z_4=0\,,\\
  z_3^5+z_4^2+z_1z_2^4=0\,.
\end{array}
%%\end{equation}
$$
This system is always of ``strict splice type'' by our comments
above.
\end{example}
\emph{Equisingular deformations} of systems of equations of strict splice
type should come from adding terms of greater or equal weight with
respect to the vertex weights to each equation.  If only greater
weight is allowed the result always is an equisingular
deformation. We speak of a \emph{higher weight deformation} and say
the resulting equations are simply of \emph{splice type}.  See
\cite{neumann-wahl10} for a fuller discussion.

The importance of splice type singularities is indicated by a
result:
\begin{theorem}{\rm\cite{neumann-wahl10}}\qua \label{th:main1} A
  system of equations of splice type defines an isolated complete
  intersection surface singularity whose link is the homology sphere
  $\Sigma$ defined by the splice diagram $\Delta$, and whose
  resolution graph is therefore the corresponding resolution diagram.

Each node $v$ of the splice diagram corresponds to an exceptional
curve $E_v$ of the resolution, and the $v$--weight of $z_i$ is its
value for the valuation given by order of vanishing on $E_v$.

  Moreover, the curve $z_i=0$ cuts out in $\Sigma$ the knot %$K_i$
  corresponding to the $i$-th leaf of $\Delta$.
\end{theorem}

One could expand the definition of
strict splice type singularities to include (for fixed $v$)
suitable linear combinations of all possible admissible monomials
associated to edges at $v$.  But, up to higher weight deformations,
this adds no generality.  Also, if we change our choice of admissible
monomials for the edges at each node, then we only change our splice
type singularities up to higher weight deformation. Thus the concept
of splice type is independent of choices of admissible monomials.
%\begin{conjecture}[Splice Type Conjecture]\label{conj:2}
%  Any Gorenstein surface singularity with integral homology sphere
%  link is a
%  complete intersection of splice type.
%\end{conjecture}
%Implicit in this conjecture is
%Conjecture \ref{conj:sg} on the necessity of the semigroup condition.

The theorem implies that the embedding dimension of a splice diagram
singularity is at most the number of leaves of the splice diagram --- but
it may be less. There are even unexpected hypersurface examples.
\begin{example}
 Let $\Delta$ be the splice diagram:
$$\splicediag{15}{20}{
\lefttag\Circ y {6pt}\lineto[drr]^(.75){q}&&&&&&&&\righttag\Circ
z {6pt}\lineto[dll]_(.75){p'}\\
&&\Circ\lineto[rr]^(.25){p''q'}^(.75){p}&&\Circ
\lineto[rr]^(.25){p'}^(.75){pq''}\lineto[dd]^(.25){p''}^(.75){p'qr}&&\Circ\\
\lefttag\Circ x {6pt}\lineto[urr]_(.75)p&&&&&&&&\righttag\Circ
w {6pt}\lineto[ull]^(.75){q'}\\
&&&&\Circ\lineto[ddl]_(.25){q''}\lineto[ddr]^(.25){p''}\\ \\
&&&\lefttag\Circ v {6pt}&&\righttag\Circ u{6pt}
}
$$
The integers $p$, $q$, $p'$, $q'$, $p''$, $q''$, $r$ are $\ge2$ and
satisfy appropriate relative primeness conditions, as well as edge
inequalities $$q'>p'q,\quad q''>p''q',\quad qr>pq''\,.$$
Associating variables $x,y,z,w,u,v$ to the leaves in clockwise order
starting from the left as shown, one may write splice equations:
\begin{align*}
  x^p+y^q&=z\\
z^{p'}+w^{q'}&=u\\
u^{p''}+v^{q''}&=x^r\\
y+w&=v
\end{align*}
These define the hypersurface singularity given by
$$((x^p+y^q)^{p'}+w^{q'})^{p''}+(y+w)^{q''}=x^r\,.$$
\end{example}

A variant of the splice diagram yields a more familiar
object. Let $\Delta$ be a splice diagram satisfying the semigroup
conditions, and choose a distinguished leaf $w'$, to form a ``rooted
diagram.'' Attach a variable $z_w$ to each leaf $w\ne w'$.  Now, for
each vertex $v$ of $\Delta$, form the same equations as before, except
that one does not consider the edge in the direction of $w'$. (One is
in general eliminating more monomials than simply setting $z_{w'}=0$
in our previous splice diagram equations.)  There is now one equation
less than there are variables. Note that the edge-weights in the
direction of $w'$ now play no role and can be discarded.  We claim
these equations generate a complete intersection curve, and this curve
is the monomial curve associated to a semigroup $\Gamma'$.  To
describe this we first briefly recall some terminology about
semigroups (see Section \ref{sec:semigroup} for more details).

The semigroups arising in this paper are always \emph{numeric
  semigroups}, that is subsemigroups $\Gamma$ of $\N=\Z_{\geq 0}$ for
which $\N-\Gamma$ is finite.  The \emph{conductor} $c(\Gamma)$ is the
smallest $c\geq 0$ so that $\gamma\geq c$ implies $\gamma\in \Gamma$.
The \emph{semigroup ring} $\C [t^{\Gamma}]$, or \emph{monomial curve}
associated to $\Gamma$, is the graded subalgebra of $\C [t]$ generated
by $t^{\gamma}, \gamma \in \Gamma$.  $\Gamma$ is called a
\emph{complete intersection semigroup} if $\C [t^{\Gamma}]$ is a
graded complete intersection.

In our situation of a splice diagram $\Delta$ satisfying the semigroup
conditions with distinguished leaf $w'$, the semigroup $\Gamma'$ is
the semigroup generated by $\ell_{w'w}$ over all leaves $w\ne w'$. We
will see in section \ref{sec:semigroup} that:
\begin{fact}
$\Gamma'$ is a complete
intersection semigroup and the modified splice equations described
above define the monomial curve $\C[t^{\Gamma'}]$.
\end{fact}

In terms of Theorem \ref{th:main1}, the significance of this curve is
that if $w'$ is the $i$-th leaf of $\Delta$ then this curve, or an
equisingular deformation of it, arises as the curve cut out by the
hyperplane $z_i=0$.

\begin{example}
  Consider the splice diagram at the beginning of Section
  \ref{sec:sd}, and let $w'$ be the lower left leaf. In the modified
  splice diagram, the weights $3$ and $11$ are removed.  Denote the
  three leaves by $w_i$, $i=1,3,4$, starting at the upper left and
  going counterclockwise; the corresponding
  variables by $z_i$; and the two nodes by $v$ and $v'$.  Then the
  equations at $v$ resp.\ $v'$ could be $z_1^2+az_3z_4=0$ and
  $z_3^5+bz_4^2=0$. The semigroup $\Gamma'$ is
  $\Gamma'=\N\langle7,4,10\rangle$ and, if we choose $a=b=-1$, the
  curve can be parametrized as $(z_1,z_3,z_4)=(t^7,t^4,t^{10})$.
\end{example}

A leaf $w'$ of a splice diagram $\Delta$ always represents a knot
in the corresponding homology sphere, and this knot is a fibered
knot (see \S11 of \cite{eisenbud-neumann}). If the homology sphere
is given as a link of a splice type singularity as above, then
this knot is the link of the curve cut out by a coordinate
hyperplane $z_i=0$ (and the fibration can be given by the usual
Milnor fibration $z_i/|z_i|$). The first Betti number of its fiber
is the \emph{Milnor number} of the knot.  We recall that even
without the semigroup condition, we have:
\begin{theorem}{\rm\cite[\S11]{eisenbud-neumann}}\qua\label{th:milnor}
The Milnor number of the above knot is $$1+\sum_{v\ne
  w'}(\delta_v-2)\ell_{vw'}\,.$$
\end{theorem}
If the link is given by splice type equations, then the theory of
curve singularities implies that this number equals the conductor of
the above semigroup $\Gamma'$, as can be confirmed by computation of
the conductor (Theorem \ref{th:cisg}).

\section{Numerical semigroups and monomial curves}
\label{sec:semigroup}
In this section we develop some results about semigroups and their
associated curves that are
needed in the proofs of Theorems \ref{thi:ends} and \ref{thi:2node} of
the Introduction.

As mentioned in Section \ref{sec:equations},
the semigroups we consider are always \emph{numeric semigroups},
that is, subsemigroups $\Gamma$ of $\N=\Z_{\geq 0}$ for which
$\N-\Gamma$ is finite.  The \emph{semigroup ring} $\C [t^{\Gamma}]$,
or \emph{monomial curve} associated to $\Gamma$, is the graded
subalgebra of $\C [t]$ generated by $t^{\gamma}, \gamma \in \Gamma$.
We briefly collect some known facts and terminology (eg,
\cite{delorme, herzog, herzog-kunz, watanabe}).

The \emph{conductor} $c(\Gamma)$ is the smallest $c\geq 0$ so that
$\gamma\geq c$ implies $\gamma\in \Gamma$.  $\Gamma$ is
\emph{symmetric} when $\gamma\in \Gamma$ if and only if
$c(\Gamma)-1-\gamma\notin \Gamma$; equivalently, $\C [t^{\Gamma}]$ is
Gorenstein (see \cite{herzog-kunz} Prop.\ 2.21).  Since $\gamma$ and
$c(\Gamma)-1-\gamma$ cannot both be in $\Gamma$, a symmetric semigroup
is maximal with given conductor.  Classically an element of $\N$ that
is not in $\Gamma$ is called a \emph{gap}. The number of gaps is
denoted $\delta(\Gamma)$; clearly
$$\delta(\Gamma)\ge c(\Gamma)/2,\quad\text{with equality if and only
  if $\Gamma$ is symmetric.}$$
$\Gamma$ is called a \emph{complete intersection semigroup} if $\C
[t^{\Gamma}]$ is a graded complete intersection.  A complete
intersection semigroup is symmetric.
$\Gamma$ is a complete intersection semigroup if and only if it
has a semigroup presentation of deficiency one (ie, with one fewer
relations than generators; see \cite{herzog}). If $\Gamma$ (complete
intersection or not) has a semigroup presentation $$\Gamma=\langle
x_1,\dots,x_n:\sum_ja_{ij}x_j=\sum_jb_{ij}x_j,i=1,\dots,r\rangle$$
with $a_{ij},b_{ij}\in\N$, then the monomial curve is presented as
$$\C[z_1,\dots,z_n]/\Bigl(\prod_jz_j^{a_{ij}}-\prod_jz_j^{b_{ij}},
i=1,\dots,r\Bigr).$$

\begin{example*}
  Relatively prime $p$ and $q$ generate a complete
intersection semigroup with conductor $(p-1)(q-1)$. This semigroup has
semigroup presentation $\langle x_1,x_2:qx_1=px_2\rangle$. Its
monomial curve $\C[t^p,t^q]$ is presented as
$\C[z_1,z_2]/(z_1^q-z_2^p)$, with the isomorphism given by $z_1\mapsto
t^p$, $z_2\mapsto t^q$.
\end{example*}

Let $(\Delta, w')$ be a finite rooted tree (tree with one vertex
singled out as ``root''), whose root vertex $w'$ is of valency 1. We
visualize it with the root vertex at the top, so ``downward'' means in
the direction away from the root.  We assume also that $\Delta$ has
positive integer weights on all edges other than the root edge and
that the weights on the downward edges at each non-root vertex are
pairwise coprime.  For example, one obtains such a tree if one picks
some leaf $w'$ of a splice diagram as root, and then forgets all ``far
weights'' of the edges of the splice diagram (from the point of view of $w'$);
equivalently, one forgets the ``near weights'' around each node.

In such a tree, the numbers $\ell_{w'v}$ for $v\ne w'$ are still defined
(product of weights on edges directly adjacent to the shortest path
from $w'$ to $v$).  We define the \emph{semigroup of $(\Delta,w')$} to
be the semigroup
$$\sg(\Delta)=\sg(\Delta,w'):=\N\langle \ell_{w'w}: w \text{ is a leaf
  of } \Delta\rangle$$
(we use the shorter $\sg(\Delta)$ if the root
vertex is clear).  Each non-root vertex of $\Delta$ cuts off a
collection of subtrees below it.  We say that $(\Delta,w')$ satisfies
the \emph{semigroup condition} if the weight on the root edge of every
such subtree is in the semigroup of the subtree.

Define an invariant $\mu(\Delta,w')$ by
$$\mu(\Delta)=\mu(\Delta,w'):=1+\sum_{v\ne w'} (\delta_v-2)\ell_{w'v}.$$
\begin{theorem}\label{th:cisg} Let $(\Delta,w')$ be a weighted rooted
  tree as above and $\Gamma=\sg(\Delta)$.  Then
  $$2\delta(\Gamma)\le \mu(\Delta),$$
  with equality if and only if
  $(\Delta,w')$ satisfies the semigroup condition, in which case
  $\Gamma=\sg(\Delta)$ is a complete intersection semigroup.  (It
  follows that the same result holds with $2\delta(\Gamma)$ replaced
  by $c(\Gamma)$.)
\end{theorem}

If $(\Delta,w')$ satisfies the semigroup condition we will describe
the complete intersection equations; these equations will be
associated to the nodes of $\Delta$. We assign a variable $z_j$ to
each leaf $w_j$ of $\Delta$. The equations will generate the kernel of
the map $\C[z_1,\dots,z_m]\to\C[t^\Gamma]$ given by $z_j\mapsto
t^{l_{w'w_j}}$.

For a node $v$ of the tree and a leaf $w_j$ below it let $l'_{vw_j}$
be the product of weights adjacent to the path from $v$ to $w_j$,
excluding weights adjacent to $v$. For each downward edge $e$ at $v$ the
semigroup condition tells us that the weight $p_e$ is a non-negative
integer linear combination $p_e=\sum_j\alpha_jl'_{vw_j}$, summed over
the leaves below $v$. We choose such an expression and denote by
$M_e=\prod_j z_j^{\alpha_j}\in\C[z_1,\dots,z_m]$ the corresponding
monomial.  Then:
\begin{scholium}
  If $(\Delta,w')$ satisfies the semigroup condition in the above
  theorem then the equations associated to node $v$ are the equations
  that equate the monomials $M_e$ for the different downward edges at
  $v$.

  If we replace each of these equations $M_e=M_{e'}$ by an equation
  $M_e=a_{ee'}M_{e'}$ with $a_{ee'}\in\C^*$ then we obtain the same
  monomial curve.
\end{scholium}
\begin{remark}
  Delorme's Proposition 9 in \cite{delorme} implies that every
  complete intersection semigroup arises as in Theorem
  \ref{th:cisg}. Already in the three-generator case the minimal tree
  defining the semigroup need not be unique.
\end{remark}
\begin{example*}
  If $\gcd(a,b)=\gcd(a,c)=\gcd(c,d)=1$ then the tree
$$\splicediag{12}{12}{&&\Circ\\
&&\Circ\lineto[u]\\
&\Circ\lineto[ur]^(.6){ad}&&\Circ\lineto[ul]_(.6){c}\\
\Circ\lineto[ur]^(.6)a&&\Circ\lineto[ul]_(.6)b\\
}
$$
satisfies the semigroup condition and leads to the complete
intersection monomial curve $$\C[z_1,z_2,z_3]/(z_1^a-z_2^b,
z_2^d-z_3^c)\cong\C[t^{bc},t^{ac},t^{ad}].$$ Exchanging $a$ with $c$ and
$b$ with $d$ gives a different tree for the same semigroup.
\end{example*}

\begin{proof}[Proof of Theorem \ref{th:cisg} and Scholium]
  The second part of the scholium is an easy induction once the rest
  is proved, replacing $z_j\mapsto t^{\ell_{w'w_j}}$ for $j>1$ by
  $z_j\mapsto \lambda_jt^{\ell_{w'w_j}}$ for suitable
  $\lambda_j\in\C^*$. So we will just prove the theorem and first part
  of the scholium.

  Let $\Delta_1,\dots,\Delta_n$ be the subtrees cut off by the bottom
  vertex $w_0$ of the root edge of $\Delta$ and let $p_i$ be the
  weight on the root edge of $\Delta_i$.  Write
  $\Gamma_i=\sg(\Delta_i,w_0)$, $P=p_1\dots p_n$ and $P_i=P/p_i$. Then
\[\Gamma =P_1\Gamma_1 +\dots+P_n\Gamma_n,\]
the semigroup consisting of all integers of the form $\sum
P_{i}\gamma_{i}$, $\gamma_{i}\in \Gamma_{i}$.
Moreover,
$$\mu(\Delta,w')= \sum_{i=1}^n\left(
    P_i\mu(\Delta_i,w_0) -1\right) + (n-1)P+1.$$
  By Lemma \ref{le:sg} below, the desired results now hold for
  $\Delta$ if they are true for each $\Delta_i$.  The proof is thus an
  induction, with the induction start being the case that $\Delta$
  consists of only a root edge and $\sg(\Delta)$ is the one-generator
  semigroup $\N$.
\end{proof}
\begin{lemma} \label{le:sg}
Suppose $\Gamma_i$ are semigroups for $i=1,\dots, n$, and
$p_1,\dots,p_n$ are pairwise coprime positive integers.  Write
$P=p_1\dots p_n$ and $P_i=P/p_i$. Let
\[\Gamma =P_1\Gamma_1 +\dots+P_n\Gamma_n.\]
Then:
\begin{enumerate}
\item\label{it:delta1} $2\delta(\Gamma)\leq \sum_{i=1}^n
  P_i(2\delta(\Gamma_i)-1) + (n-1)P +1$.
\item\label{it:delta2} If equality holds in {\rm(\ref{it:delta1})} then
  $p_i\in \Gamma_i$ for $i=1,\dots,n$.
\item\label{it:sg1} $c(\Gamma)\leq \sum_{i=1}^n P_i\left(c(\Gamma_i)
    -1\right) + (n-1)P+1.$
\item\label{it:sg2} If $p_i\in \Gamma_i$ for $i=1,\dots,n$ then
  equality holds in {\rm(\ref{it:sg1})}.
\item\label{it:sg2.2} If each $\Gamma_i$ is symmetric then the three
  statements are equivalent: equality in {\rm(\ref{it:delta1})};
  equality in {\rm(\ref{it:sg1})}; $p_i\in\Gamma_i$ for $i=1,\dots,n$.
\item\label{it:sg3} Assuming $p_i\in \Gamma_i$ for each $i$, then
  $\Gamma$ is symmetric resp.\ a complete intersection if and only if
  each $\Gamma_i$ is symmetric resp.\ complete intersection.
\item\label{it:sg2.5} If $p_i\in \Gamma_i$ for each $i$ then one
  obtains a presentation for $\Gamma$ by adjoining to the disjoint
  union of presentations for the $\Gamma_i$ the $n-1$ relations
  $w_1=\dots=w_n$, where $w_i$ is an expression for $p_i$ in the
  presentation of\/ $\Gamma_i$.
\end{enumerate}
\end{lemma}
\begin{proof} We shall prove the case $n=2$. The case of general $n$
  follows from this case by an easy induction.

To prove (\ref{it:delta1}) we count gaps in $\Gamma$.
  A gap $\gamma$ of $\Gamma=p_2\Gamma_1+p_1\Gamma_2$ is either
  \begin{enumerate}
  \item[(i)] one of
  the $(p_1-1)(p_2-1)/2$ gaps of $p_2\N+p_1\N$,
  \end{enumerate}
or it is of the form
  $\gamma=p_2\alpha+p_1\beta$ for some $\alpha,\beta\in\N$. In this
  case we will see that either:
\begin{enumerate}
\item[(ii)] $\beta$ is the smallest %number satisfying
  $\beta\in\Gamma_2$ in its congruence class mod $p_2$, and
  $\alpha\notin\Gamma_1$, or
\item[(iii)] $0\le\alpha<p_1$, and $\beta\notin\Gamma_2$.
\end{enumerate}
Indeed, if we can express $\gamma$ in the form
$\gamma=p_2\alpha+p_1\beta$ with $\alpha,\beta\in\N$, then we can do
so with $0\le\alpha<p_1$. If this expression does not satisfy
condition (iii) then $\beta\in\Gamma_2$. In this case decrease $\beta$
by some multiple of $p_2$ (maybe zero) to make it the smallest
$\beta\in\Gamma_2$ in its congruence class mod $p_2$, and
simultaneously increase $\alpha$ by the same multiple of $p_1$ to keep
$\gamma=p_2\alpha+p_1\beta$. Since $\gamma$ is a gap of
$p_2\Gamma_1+p_1\Gamma_2$, we must have $\alpha\notin\Gamma_1$, so the
expression now satisfies condition (ii).

Now there are exactly $p_2\delta(\Gamma_1)$ pairs $(\beta,\alpha)$
satisfying condition (ii) and $p_1\delta(\Gamma_2)$ pairs satisfying
condition (iii), so there are at most
$(p_1-1)(p_2-1)/2+p_2\delta(\Gamma_1)+p_1\delta(\Gamma_2)$ gaps of
$\Gamma=p_2\Gamma_1+p_1\Gamma_2$. This number can be written
$\frac12\bigl(p_2(2\delta(\Gamma_1)-1)+
p_1(2\delta(\Gamma_2)-1)+p_1p_2+1\bigr)$, so part (\ref{it:delta1}) is
proven.

This proof shows that we have equality in part (\ref{it:delta1}) if
and only if every element $\gamma=p_2\alpha+p_1\beta$ satisfying
condition (ii) or (iii) is a gap of $\Gamma$ and there is no overlap
between cases (ii) and (iii).  Suppose now $p_1\notin\Gamma_1$.  Then if
every $p_2\alpha+p_1\beta$ satisfying (ii) is a gap of $\Gamma$, there
is an overlap: $(\alpha,\beta)=(p_1,0)$ in condition (ii) shows that
$p_1p_2$ is a gap of $\Gamma$, whence $p_2\notin\Gamma_2$, so $p_1p_2$
also has an expression with $(\alpha,\beta)=(0,p_2)$ satisfying
condition (iii).  Thus $p_1\notin\Gamma_1$ implies inequality in part
(\ref{it:delta1}).  Similarly for $p_2\notin\Gamma_2$, so part
(\ref{it:delta2}) is proved.

For statement (\ref{it:sg1}), we show that $i\geq 0$ added to the
right hand side of the inequality of part (\ref{it:sg1}) gives an
element of $\Gamma$.  The sum of the last two terms of $$p_2
c(\Gamma_1)+p_1 c(\Gamma_2) +(p_1-1)(p_2-1) +i$$
is in the semigroup
generated by $p_1$ and $p_2$, say $p_1\alpha +p_2\beta$; so the whole
expression equals
\[p_2(c(\Gamma_1)+\beta)+p_1(c(\Gamma_2)+\alpha),\]
which by definition of conductors is clearly in $\Gamma$.

For statement (\ref{it:sg2}), suppose $p_1\in \Gamma_1$ and $p_2\in
\Gamma_2$, but
   \[p_2c(\Gamma_1) + p_1c(\Gamma_2)+(p_1-1)(p_2-1)-1
   =p_2\lambda +p_1\pi,\text{ for some } \lambda\in \Gamma_1, \pi\in
   \Gamma_2.\] Modulo $p_1$ this equation says
   $c(\Gamma_1)-1\equiv\lambda$, so
   \[c(\Gamma_1)-1=\lambda+p_1t,\text{ for some integer }t.\]
Inserting this in the previous equation gives
   \[c(\Gamma_2)-1=\pi +p_2(-1-t).\]
   Since one of $t$ and $-1-t$ is $\ge0$ and $\lambda, p_1 \in \Gamma_1$
   and $\pi,p_2 \in \Gamma_2,$ one gets either $c(\Gamma_1)-1\in\Gamma_1$ or
   $c(\Gamma_2)-1\in \Gamma_2$, a contradiction.

   Part (\ref{it:sg2.2}) is now immediate: (\ref{it:delta2}) and
   (\ref{it:sg2}) show
   $$\text{(equality in (\ref{it:delta1})) }\quad\Rightarrow\quad
   (p_1\in\Gamma_1\text{ and }p_2\in\Gamma_2) \quad\Rightarrow\quad
   \text{(equality in (\ref{it:sg1})), }$$
   and if the $\Gamma_i$ are
   symmetric then $c(\Gamma)\le2\delta(\Gamma)$ and
   $c(\Gamma_i)=2\delta(\Gamma_i)$, so equality in (\ref{it:sg1})
   implies equality in (\ref{it:delta1}).

   Part (\ref{it:sg3}) is proved in \cite{delorme}. (In this paper we
   use only that $\Gamma$ is a complete intersection if both
   $\Gamma_1$ and $\Gamma_2$ are; this follows from part
   (\ref{it:sg2.5}).)

   For part (\ref{it:sg2.5}), let $\Gamma_1=\langle
   x_1,\dots,x_n:s_1,\dots,s_{k}\rangle$ and $\Gamma_2=\langle
   y_1,\dots,y_m:r_1,\dots,r_{\ell}\rangle$ be commutative semigroup
   presentations of\/ $\Gamma_1$ and $\Gamma_2$;  let
   $p_1=v(x_1,\dots,x_n)$ and $p_2=w(y_1,\dots,y_m)$ be expressions
   for $p_1$ and $p_2$ in these semigroups. Suppose
   $p_2\gamma_1+p_1\gamma_2=p_2\gamma_1'+p_1\gamma_2'$ equates two
   elements of $\Gamma=p_2\Gamma_1+p_1\Gamma_2$, with
   $\gamma_1,\gamma_1'\in \Gamma_1$ and
   $\gamma_2,\gamma_2'\in\Gamma_2$. Let $\gamma_1=g_1(x_1,\dots,x_n)$
   be an expression for $\gamma_1\in\Gamma_1$ in terms of the
   generators (and hence for $p_2\gamma_1$ in $p_2\Gamma_1)$, and
   similarly $\gamma_1'=g_1'(x_1,\dots,x_n),
   \gamma_2=g_2(y_1,\dots,y_m), \gamma_2'=g_2'(y_1,\dots,y_m)$.  Then
   the relation to be verified in $\Gamma=p_2\Gamma_1+p_1\Gamma_2$ is
   $g_1+g_2=g_1'+g_2'$ (abbreviating $g_1(x_1,\dots,x_n)=g_1$ etc.),
   and we must show this follows from the relations of $\Gamma_1$ and
   $\Gamma_2$ and the additional relation $v=w$.

   With no loss of generality $\gamma_1\ge\gamma_1'$ in $\N$.  Then,
   working in $\N$, we have
   $p_2(\gamma_1-\gamma_1')=p_1(\gamma_2'-\gamma_2)$, so
   $\gamma_1-\gamma_1'=sp_1$ and $\gamma_2'-\gamma_2=sp_2$ for some
   $s$ in $\N$. In particular, the equations $g_1=sv+g_1'$ and
   $g_2'=sw+g_2$ hold in $\Gamma_1$ and $\Gamma_2$, so they must
   follow from the relations of these semigroups. Thus, using the
   additional relation $v=w$, we deduce
   $g_1+g_2=sv+g_1'+g_2=sw+g_1'+g_2=g_1'+g_2'$, as desired.
  \end{proof}

\subsection{Normal form monomials}\label{subsec:normalform}
The material of this subsection will be needed in Section
\ref{sec:genus} for the proof of Theorem
\ref{thi:2node}.

Suppose now that $(\Delta, w')$
satisfies the semigroup condition and put $\Gamma=\sg(\Delta)$. We
wish to describe a monomial basis for the corresponding complete
intersection curve $\C[z_1,\dots,z_m]/($relations$)$. That is, we want
``normal form'' monomials in $z_1,\dots,z_m$ so that each $t^\gamma$
with $\gamma\in \Gamma$ is the image of exactly one monomial under the
map $\C[z_1,\dots,z_m] \to \C[t^\Gamma]$ given by $z_j\mapsto
t^{\ell_{w'w_j}}$. We will do this by systematically trying to
eliminate variables with small index.

We assume that the tree $\Delta$ is drawn so that the indices
$i=1,\dots,m$ of the leaves increase from left to right.  For any node
$v$ and outward edge $e$ at $v$ let $\Delta_{ve}$ be the subtree below
$v$ with root vertex $v$ and root edge $e$.

If $M$ is a monomial, let $M_{ve}$ be the submonomial of $M$
determined by the variables corresponding to leaves of $\Delta_{ve}$.
This monomial represents $t^\alpha\in\C[t^{\sg(\Delta_{ve})}]$ for
some $\alpha$.  We will say $M$ is in \emph{normal form} if for every
$v$ and $e$ as above so that $e$ is not the rightmost edge at $v$,
$\alpha-p_e\notin\sg(\Delta_{ve})$.

If $M$ is not in normal form at some $(v,e)$ then we could replace
$M_{ve}$ in $M$ by $M'_{ve}M_{e'}$ where $e'$ is the rightmost edge at
$v$, $M'_{ve}$ is a monomial representing
$t^{\alpha-p_e}\in\C[t^{\sg(\Delta_{ve})}]$ and $M_{e'}$ is a monomial
representing $t^{p_{e'}}\in\C[t^{\sg(\Delta_{ve'})}]$.  Since
$t^{p_e}\in\C[t^{\sg(\Delta_{ve})}]$ and
$t^{p_{e'}}\in\C[t^{\sg(\Delta_{ve'})}]$ become equal in
$\C[t^{\sg(\Delta)}]$, this does not change the value of $M$. It is
easy to see this process must eventually stop. A simple induction
shows that it yields a unique normal form for $M$.  Normal form
monomials thus provide the desired monomial basis of
$\C[z_1,\dots,z_m]/($relations$)$.

The following example will be important in Section \ref{sec:genus}.
\begin{example}\label{ex:monomialbasis} Let
$$\Delta=\splicediag{17}{14}{
&\Circ\lineto[d]\\
&\Circ\lineto[dl]_(.4){p_1}\lineto[dr]^(.4){p_n}\\
\Circ&\dots&\Circ}
$$
so $\Gamma$ is the semigroup generated by the $P_{i}=P/p_i$.
The monomial curve \[(t^{P_{1}},t^{P_{2}},\cdots ,t^{P_{n}})\]
    is the complete intersection curve singularity defined by the
    equations \[z_{i}^{p_{i}}-z_{n}^{p_{n}}=0,\quad  i=1,\cdots ,n-1.\]
The conductor $c(\Gamma)$ is \[P\left(n-1-\sum (1/p_{i})\right)+1.\]

The monomial basis described above is
$$\{z_1^{\alpha_1}\dots z_n^{\alpha_n}:
\alpha_{i}< p_i\text{ for all }i=1,\dots,n-1\}.$$

More generally, applied to a  tree of the form
$$\Delta=\splicediag{10}{8}{
\Circ\lineto[ddrr]\\
\\
&&\Circ\lineto[ddll]_(.4){p_1}\lineto[dd]^(.6){p_{k_1}}
\lineto[dddrrrr]^(.5){q_1}\\ \\
\Circ&\dots&\Circ\\
&&&&&&\Circ
\lineto[ddll]_(.4){p_{k_1+1}}\lineto[dd]\dashto[dddddrrrrrr]\\ \\
& &&& \Circ&\dots &\Circ\\ \\ \\
&&&&& & &&& & & &
\Circ\lineto[ddl]_(.4){p_{k_r+1}}
\lineto[ddr]^(.4){p_n}\\ \\
&&&&&&& &&&& \Circ &\dots &\Circ
}
$$
which satisfies the semigroup condition, the above procedure will
again give the monomial basis
$$\{z_1^{\alpha_1}\dots z_n^{\alpha_n}: \alpha_{i}< p_i\text{ for all
  }i=1,\dots,n-1\}.$$
(However, with a different ordering of the
variables the monomial basis for this example can be considerably more
complicated.)
\end{example}

\section{The semigroup condition}\label{sec:semigroup condition}
Let $(X,o)$ be a normal surface singularity whose link $\Sigma$ is an
integral homology sphere.  Each leaf of the splice (or resolution)
diagram gives a knot in $\Sigma$, unique up to isotopy.  A key point
in the proof in \cite{neumann-wahl10}, that splice diagram equations
give integral homology sphere links, is to show that the variable
$z_{i}$ associated to a leaf cuts out the corresponding knot in
$\Sigma$.  In other words, the curve $C_{i}$ given by $z_{i}=0$ is
irreducible, and its proper transform $D_{i}$ on the minimal good
resolution is smooth and intersects transversely the exceptional curve
corresponding to the leaf of the splice diagram.  We show that the
existence of such functions implies the semigroup condition on the
splice diagram.

\begin{theorem} \label{th:ends}
  Let $(X,o)$ be
  a normal surface singularity whose link $\Sigma$ is an integral
  homology sphere.  Assume that for each of the $t$ leaves $w_i$ of
  the splice diagram $\Delta$ of\/ $\Sigma$, there is a function
  $z_{i}$ inducing the end knot as above.  Then
  \begin{enumerate}
   \item $\Delta$ satisfies the semigroup condition
   \item $X$ is a complete intersection of embedding dimension $\leq t$
   \item $z_{1},\cdots,z_{t}$ generate the maximal ideal of the local
     ring of $X$ at $o$, and $X$ is a complete intersection of splice
     type with respect to these generators.% (up to higher weight deformation).
  \end{enumerate}
\end{theorem}
\begin{proof}  Let
  $(Y,E)\rightarrow (X,o)$ be the minimal good resolution, $z=z_{1}$ a
  function as above, $C\subset X$ the irreducible Cartier divisor
  defined by $z=0$, $D\subset Y$ its proper transform, and
  $E_{1}\subset Y$ the exceptional curve (which intersects $D$ in one
  point) corresponding to the leaf of the splice diagram.

  Let $V$ be the value semigroup of $C$. The orders of vanishing of
  the functions $z_{2},\cdots,z_{t}$ at $D\cap E_{1}$ generate a
  subsemigroup $\Gamma\subset V$ which we can compute from $\Delta$ as
  follows.  For each exceptional curve $E_i$, let $a_{ij}$ be the
  order of vanishing of $z_j$ on $E_i$, so, as a divisor,
  $(z_j)=\sum_ia_{ij}E_i+D_j$. The equations $z_j^{-1}(0)\cdot
  E_k=0$ imply that $a_{ij}$ is the $ij$--entry of the matrix
  $(-E_i\cdot E_j)^{-1}$; so $a_{ij}=\ell_{ij}$ (see Theorem
  \ref{th:props}).
  Thus $\Gamma$ is
  the semigroup generated by $\ell_{1j}$, $j\ge2$.

  Theorem \ref{th:cisg} implies $2\delta(\Gamma)\le\mu(\Delta, w_1)$,
  where $\mu(\Delta,w_1)$ is described there and $\delta(\Gamma)$
  denotes the number of gaps of $\Gamma$. But, by Theorem
  \ref{th:milnor}, $\mu(\Delta,w_1)$ is also equal to the
  $\mu$--invariant $\mu(C)$ of the curve $C$.  Now $\mu(C)=2\delta(V)$
  (since we do not know a priori that the curve is Gorenstein, we must
  appeal to Buchweitz and Greuel \cite{buchweitz-greuel} for this).
  Since the inclusion $\Gamma\subset V$ implies
  $\delta(V)\le\delta(\Gamma)$, we conclude that
  $2\delta(V)=2\delta(\Gamma)=\mu(\Delta,w_1)$. Thus $\Gamma=V$, and,
  by Theorem \ref{th:cisg} again, $\Gamma=V$ is a complete
  intersection semigroup.  This implies that $C$ is a positive weight
  deformation of the monomial curve $\C [t^{\gamma}: \gamma \in
  \Gamma]$ (eg, Teissier's appendix to \cite{zariski} or
  \cite{teissier}) and in particular is itself a complete intersection
  (with maximal ideal generated by the images of
  $z_{2},\cdots,z_{t}$).  It follows that $(X,o)$ is a complete
  intersection (with maximal ideal generated by $z_1,\dots,z_n$).
  Finally, repeating the argument at every leaf gives all the
  semigroup conditions.

  It remains to show that, using the functions $z_{1},\cdots,z_{t}$
  above, we can find splice equations for the singularity.
  This will proceed as follows: for each node $v$ of
  valency $\delta=\delta_v$, we will write down appropriate monomials
  in the $z_{i}$ which have the same weight at the node (ie, order
  of vanishing along the corresponding exceptional curve), and
  conclude there are $\delta -2$ independent linear dependence
  relations among these monomials, mod higher weight terms.

  Let $E_v$ be the exceptional curve corresponding to the node $v$ and
  let $E_1,\dots,E_\delta$ be the exceptional curves with intersect
  $E_v$, corresponding to edges $e_1,\dots,e_\delta$ at $v$.
  Choose a monomial $M_i$ of weight $d_v$ associated to each edge
  $e_i$ at $v$ (their existence is guaranteed by the semigroup
  condition). On the exceptional curve $E_v$ these monomials all
  vanish to order $d_v$. If we go to an adjacent node $v'$ of the
  splice diagram, as in
  $$\splicediag{8}{30}{
    &&&&\\
    \Vdots&\overtag\Circ v {8pt}\lineto[ul]_(.5){p_1}
    \lineto[dl]^(.5){p_{\delta-1}}
    \lineto[rr]^(.25){p_\delta}^(.75){q_{\delta'}}&&
    \overtag\Circ{v'}{8pt} \lineto[ur]^(.5){q_1}
    \lineto[dr]_(.5){q_{\delta'-1}}&\Vdots\\
    &&&& }$$
  then the order of vanishing of $M_i$ is $p_1\dots
  p_{\delta-1}q_1\dots q_{\delta'-1}$ for $i\ne\delta$ and the order
  of vanishing of $M_\delta$ is $p_\delta q_{\delta'}$. In particular,
  at the exceptional curve corresponding to $v'$, $M_{\delta}$
  vanishes to order $D$ more than the other $M_i$'s, where $D$ is the
  edge determinant of edge $e_\delta$. Now in the (unreduced) maximal
  splice diagram (see the Appendix; Section \ref{sec:splicing}) we
  have a node for every exceptional curve and all edge determinants
  are $1$.  Thus we see that on each exceptional curve $E_i$ that
  intersects $E_v$, the $M_j$ with $j\ne i$ vanish to a common order
  and $M_i$ vanishes to one higher order.  Thus, if we fix one of the
  neighboring exceptional curves, say $E_\delta$, then each ratio
  $M_{i}/M_{\delta}$ for $i\ne \delta$ gives a function on $E_v$ that
  has a pole of order 1 at the point $E_v\cap E_\delta$, a simple zero
  at the point of intersection $E_v\cap E_i$, and no other poles or
  zeros. It follows that there are $\delta-2$ linearly independent
  relations among the $M_i$ up to higher order at $E_v$, as desired.

  This gives us a collection of higher weight perturbations of
  equations of strict splice type and they are the complete intersection
  description of $(X,o)$  since they give the appropriate complete
  intersection curves when intersected with $z_j=0$.
\end{proof}

It is a Riemann-Roch problem to determine if a
%Gorenstein
singularity
with homology sphere link has functions $z$ with the properties
described above.  However, it is not even known if
there is any function at all giving an \emph{irreducible} divisor on
$X$; this is certainly not the case for a general hypersurface
singularity \cite{laufer}.

We give an application of the above theorem. We will show in
Section \ref{sec:plane curves} that if a surface singularity of the
form $z^n=g(x,y)$ has homology sphere link, then there is a splice
type singularity with the same topology (and this singularity is
analytically equivalent to one given by an equation of the form
$z^n=f(x,y)$).  This leaves open the question whether the original
singularity $z^n=g(x,y)$ is an equisingular deformation of the strict splice
type singularity $z^n=f(x,y)$ and is hence of splice type.
\begin{corollary}\label{cor:znfxy}
  Any surface singularity with homology sphere link given by an
  equation $z^n=g(x,y)$ is a splice type singularity.
\end{corollary}
\begin{proof}
  We just sketch the proof. If the splice diagram for the plane curve
  $g(x,y)=0$ is
$$
\splicediag{16}{24}{
    \Circ\lineto[r]^(.75){p_1}&\Circ\lineto[r]^(.25){1}^(.75){p_2}
    \lineto[d]^(.25){q_1}&\Circ\lineto[d]^(.25){q_2}\lineto[r]^(.25)1&
    \dotto[r]&\lineto[r]^(.75){p_k}&\Circ\lineto[d]^(.25){q_k}\ar[r]^(.25)1&\\
    &\Circ&\Circ&&&\Circ}
  $$
  then it is known (see, eg, \cite{teissier}) that curves
  corresponding to ends of this splice diagram are cut out by
  polynomials (namely certain ``approximate roots'' $g_i(x,y)$ of
  $g(x,y)$). It is easy to check that the functions $g_i(x,y)$ then
  cut out curves in the surface $z^n=g(x,y)$ corresponding to the ends
  of its splice diagram
$$\splicediag{16}{24}{
    \Circ\lineto[r]^(.75){p_1}&\Circ\lineto[r]^(.25){n}^(.75){p_2}
    \lineto[d]^(.25){q_1}&\Circ\lineto[d]^(.25){q_2}\lineto[r]^(.25)n&
    \dotto[r]&\lineto[r]^(.75){p_k}&\Circ\lineto[d]^(.25){q_k}\lineto
[r]^(.25)n&\Circ\\
    &\Circ&\Circ&&&\Circ}
$$
so Theorem \ref{th:ends} applies.
\end{proof}

\section{Geometric genus and Theorem \ref{thi:2node}}\label{sec:genus}

In this section we will prove Theorem \ref{thi:2node}, that the Casson
Invariant Conjecture holds for a splice type singularity when the
nodes of the splice diagram are in a line. We will do this by
computing geometric genus $p_g$, to prove Version 2 of the Casson
Invariant Conjecture in the Introduction. The equivalence of the two
versions of the Casson Invariant Conjecture will be proved in Theorem
\ref{th:equiv}.

Let $(X,o)$ be a germ of a normal surface singularity, with analytic
local ring $\mathcal O$.  Consider a good resolution
$\pi:(Y,E)\rightarrow (X,o)$, ie, the exceptional fiber $E=\bigcup
E_{i}$ is a union of smooth curves intersecting transversely, no three
through a point.  By local duality, one may compute the geometric
genus in two ways:
\[p_{g}(X)=\text{dim}\ H^{1}(\mathcal O_{Y})=\text{dim}\
H^{0}(U,K_{U})/H^{0}(Y,K_{Y}),\] where $U=X-\{o\}=Y-E$, and $K$
denotes canonical line bundle (or its sheaf of sections).

If $(X,o)$ is Gorenstein, let $\omega$ be a nowhere-0 holomorphic
two-form on $U$.  Define the \emph{canonical ideal} $J$ of $\mathcal
O$ by \[J=\{f\in \mathcal O : f\omega\ \text{is regular on}\ Y\}.\]
Then clearly \[p_{g}(X)= \text{dim}\ \mathcal O /J.\] Let
$E_{\alpha},\alpha=1,\ldots,t$ be those exceptional curves which
either have positive genus, or intersect at least three other curves.
Let $G$ be the union of the remaining curves (the ``strings'' in the
resolution).  The blowing-down $Y\rightarrow Y'$ of $G$ gives a space
with only cyclic quotient singularities (if $Y$ is the minimal good
resolution then $Y'$ is the ``log-canonical resolution''); since these
singularities are rational, regular forms in a punctured neighborhood
automatically extend regularly on a resolution.  Therefore, $f\in J$
if and only if $f\omega$ extends regularly over the $t$ particular
curves $E_{\alpha}$.  Let $\nu_{\alpha}$ be the valuation on $\mathcal
O$ given by order of vanishing along $E_{\alpha}$, and let
$k_{\alpha}-1$ denote the order of the pole of $\omega$ along that
curve.  We conclude that
\[J=\{f\in \mathcal O : \nu_{\alpha}(f)\geq k_{\alpha}-1, \alpha
=1,2,\ldots,t\}.\]

In our case we can improve $k_\alpha-1$ to $k_\alpha$.

\begin{proposition} Let $(X,o)$ be the germ of a Gorenstein surface
  singularity, whose link is a rational homology sphere.  Let
  $(Y,E)\rightarrow (X,o)$ be the minimal good resolution, and let
  $E_{1},\cdots ,E_{t}$ be the exceptional curves of valency $\geq 3$.
  Let $k_{\alpha}$ be the coefficient of $E_{\alpha}$ in the divisor
  $-(K+E)$, and $\nu_{\alpha}$ the corresponding valuation of the local
  ring $\mathcal{O}$ of $X$.  Then the geometric genus of $X$ is the
  colength of the ideal \[J=\{f\in \mathcal O : \nu_{\alpha}(f)\geq
  k_{\alpha},\ \alpha=1,\ldots ,t\}.\]
\end{proposition}

\begin{proof} By the preceding discussion,
the statement to be proved is
$$H^{0}(Y, K_{Y})=H^{0}(Y-G, K_{Y}+E).$$
We will do this in two steps:
$$H^{0}(Y, K_{Y})=H^{0}(Y, K_{Y}+E)=H^{0}(Y-G, K_{Y}+E).$$
Since the link of $X$ is a $\Q$--homology-sphere, the exceptional curve
$E$ is the transverse union of smooth rational curves $E_i$, no three
through a point, with contractible dual graph. It follows that
$h^1(\mathcal O_E)=0$. (Proof: write $E=E_1+F$, where $E_1$ is a
component of $E$ that meets the rest $F$ of $E$ in a single point; the
surjection $\mathcal O_E\to \mathcal O_F$ has kernel $\mathcal
O_{E_1}(-F)=\mathcal O(-1)$, so the claim follows by induction on the
number of components of $E$.)

Denote $K_Y\otimes\mathcal O_E(E)$ by $K_E$ (called the dualizing
sheaf in \cite{barth-peters-vandeven}, Section II.1).  Serre duality
implies that, for any line bundle $L$ on $E$, $H^1(E,L)$ is dual to
$H^0(E,L^*\otimes K_E)$ (eg, \cite{barth-peters-vandeven}, Theorem
II(6.1)).  Taking $L$ trivial we see $h^{0}(K_{E})=0$.  The adjunction
sequence $0\to K_Y\to K_Y+E\to K_E\to 0$ (called ``residue sequence''
in \cite{barth-peters-vandeven}, Section II.1) now gives $ 0\to
H^0(K_Y)\to H^0(K_Y+E)\to H^0(K_E)=0$, proving the first equality
$H^{0}(K_{Y})=H^{0}(K_{Y}+E)$.

The second equality $H^0(Y, K_Y+E)=H^0(Y-G, K_Y+E)$ holds generally,
without the condition on the link.  In fact, if $G$ is any union of
components of $E$ and $L$ any divisor supported on $E$, then it is easy
to see that $H^0(Y, L)=H^0(Y-G, L)$ so long as $L\cdot G_i\le 0$ for
each component of $G$ (for a stronger statement see \cite{giraud}), so
we must just show that that $(K+E) \cdot G_i\le0$ for all $i$.  But
$G_i$ is a smooth rational curve, so $(K+E)\cdot G_i$ equals $-2$ plus
the number of intersections of $G_i$ with the other curves of $E$.
This result is $-1$ if $G_i$ is an end curve of the graph, or $0$
otherwise.  In either case, the condition is fulfilled, and our result
follows.
\end{proof}

While the $k_{i}$ are determined from the resolution graph (see
Proposition \ref{prop:K}), in some cases they can be computed directly
from the equations defining $\mathcal O$.

\begin{proposition} \label{prop:k}Let
    \[\C[z_{1},\ldots,z_{s}]/(f_{1},\cdots,f_{s-2})\]
    define an isolated complete intersection surface singularity at
    the origin.  For an exceptional curve $E_{1}$ in a resolution,
    with valuation $\nu=\nu_{1}$, consider the filtration defined by
    $I_{n}=\{f : \nu(f)\geq n\}$.  Assume that the associated graded of
    this filtration is a complete intersection integral domain, with
    the $z_{i}$ inducing homogeneous generators, and defined by the
    $\nu$--leading forms $\overline{f_{j}}, j=1,\cdots,s-2$. Then the
    invariant $k_{1}$ is computed as
\[k_{1}=\sum_{j=1}^{s-2} \nu(\overline{f_{j}}) \ -\sum_{i=1}^{s} \nu(z_{i}).\]
\end{proposition}

\begin{proof}
  We may interpret \[\omega =dz_{1}\wedge \cdots \wedge
  dz_{s}/df_{1}\wedge \cdots \wedge df_{s-2}.\] On the associated
  graded, this gives a two-form of total weight
\[\Sigma \nu(z_{i}) -\Sigma \nu(\overline{f_{j}}) .\]
In terms of local coordinates in a neighborhood of a general point of
$E_{1}$, one finds the order of the pole of $\omega$ is one more than
the weight, as desired.
\end{proof}

If our singularity is a complete intersection of splice type and $E_1$
corresponds to a node $v$ of the splice diagram, then, in the
terminology of the preceding section, $\nu(z_i)$ is the $v$--weight of
$z_i$, so $\nu(z_i)$ is the product of splice diagram weights adjacent
to the path from node $v$ to leaf $i$.  It is easy to see that the
formula of the above proposition is then equivalent to that of
Proposition \ref{prop:K}.

\begin{example} The last two propositions give a well-known
  result for a weighted homogeneous complete intersection: the
  geometric genus is the sum of the dimensions of the graded pieces of
  weight less than or equal to $k_{1}$ above.  In particular, let
  $V(p_{1},\ldots,p_{n})$ (with $p_{i}$ pairwise relatively prime) be
  a \BCI{}, defined by
%\[z_{1}^{p_{1}}+a_{1}z_{n-1}^{p_{n-1}}+b_{1}z_{n}^{p_{n}}=0\]
%\[z_{2}^{p_{2}}+a_{2}z_{n-1}^{p_{n-1}}+b_{2}z_{n}^{p_{n}}=0\]
%\[ \cdots\]
\[z_{i}^{p_{i}}+a_{i}z_{n-1}^{p_{n-1}}+b_{i}z_{n}^{p_{n}}=0,\quad i=1\dots,n-2\,.\]
Let $P=p_{1}\cdots p_{n}, P_{i}=P/p_{i}$.  Then $$k_{1}=(n-2)P-\Sigma
P_{i}=P(n-2-\Sigma (1/p_{i})).$$
Using the monomial basis $z_{1}^{i_{1}}\cdots
z_{n}^{i_{n}}$ with $i_{k}<p_{k},\ k=1,\cdots,n-2$, one computes
\begin{align*}
  p_{g}(V(p_{1},\ldots,p_{n}))=
\#\{(i_{1},i_{2}, \cdots,i_{n}&)\in (\Z_{\geq 0})^{n}:
\sum_{k=1}^{n}
(i_{k}+1)/p_{k}<n-2;\\
&i_{k}<p_{k},\ k=1,\cdots,n-2\}.
\end{align*}
\end{example}

To extend this calculation to singularities corresponding to more
complicated splice diagrams, we need \emph{one} monomial basis which
works for \emph{every} filtration defined by a node of the splice
diagram.

Suppose we have a complete intersection $(X,p)$ of splice type
corresponding to a splice diagram $\Delta$. For convenience of
notation we will assume equations of strict splice type (no higher
order terms); the identical proofs will handle the general case.  Let
$\nu$ be the valuation associated to the node $v$ of $\Delta$ (see
Theorem \ref{th:main1}).  Let the edges around $v$ be $e_1, \dots,
e_n$ with weights $d_{ve_i}=p_i$, $i=1,\dots,n$ at $v$.  For each node
$v'$ of $\Delta$ the equations have the form $\sum a_{e'} M_{v'e'}=0$,
sum over the edges $e'$ at $v'$, where $M_{v'e'}$ is an admissible
monomial at $v'$ and $a_{e'}\in\C$. If $v'\ne v$ and ${e'}$ is the
edge on the path from $v'$ to $v$ we will call $M_{v'e'}$ a \emph{near
  monomial} at $v'$ for $v$.  Thus there is one near monomial for $v$
associated to each node other than $v$.
\begin{theorem}\label{th:assocgraded}
  The associated graded ring $R$ of\/ $(X,p)$ with respect to the
  filtration associated to $\nu$ is a reduced and irreducible complete
  intersection, defined by the same equations as $(X,p)$ but with the
  coefficients of all near monomials for $v$ set to zero (so only the
  equations associated to the node $v$ remain unchanged). Its
  normalization is the \BCI{} $V(p_{1},\cdots,p_{n})$.
\end{theorem}
We will need a specific basis of the graded ring $R$.
\begin{proposition}\label{prop:assocgraded}
  Choose an edge $e$ at $v$ and picture the edge $e$ as horizontal,
  with $v$ on the left.  Cut the edge $e$ at its midpoint and use this
  midpoint as root of the resulting trees $\Delta_e^L$ on the left and
  $\Delta_e^R$ on the right (so these are rooted trees and
  $\Delta_e^L$ contains $v$). Let $\mathcal M_e^L$ and $\mathcal
  M_e^R$ be monomial bases for the monomial curves of $\Delta_e^L$ and
  $\Delta_e^R$, constructed as in subsection \ref{subsec:normalform}.
  Then the set of monomials $\mathcal M_e^L\mathcal
  M_e^R=\{M_1M_2:M_1\in \mathcal M_e^L, M_2\in \mathcal M_e^R\}$ forms
  a $\C$--basis of the associated graded ring $R$.
%The normalization   is the \BCI{} $V(p_{1},\cdots,p_{n})$.
  The integer $k_{v}$ is given by
  $$k_{v}=d_{ve}(C^L-1)+\frac{d_v}{d_{ve}}(C^R-1)\quad\Bigl(=
    p_j(C^L-1)+P_j(C^R-1)\,\text{ if }e=e_j\Bigr)$$
    where $C^L$ and $C^R$ are the conductors of the semigroups
    $\sg(\Delta_e^L)$ and $\sg(\Delta_e^R)$.
\end{proposition}
\begin{proof}[Proof of Theorem \ref{th:assocgraded} and Proposition \ref{prop:assocgraded}]
  Let $\nu'$ be the valuation associated to $v'$.  If $z_w$ is the
  variable associated to a leaf $w$ then one checks easily that
  $$\frac{\nu(z_w)}{\nu'(z_w)}=\frac{\ell_{v'v}}{d_{v'}}D_{e'_1}\dots
    D_{e'_k}$$
where:
\begin{itemize}
\item $\ell_{v'v}$ is, as usual, the product of weights adjacent to
  the path from $v'$ to $v$;
\item $d_{v'}$ is the product of edge weights at $v'$;
\item $e'_1,\dots,e'_k$ are the edges that are on the path from $v'$
  to $v$ but not on the path from $w$ to $v$;
\item   for any edge $e$, $D_e$ is the product of the edge
  weights on $e$ divided by the product of edge weights directly
  adjacent to $e$ (so $D_e>1$ by the edge determinant condition).
\end{itemize}
Thus $\nu(z_w)/\nu'(z_w)$ takes its minimum value (namely
${\ell_{v'v}}/{d_{v'}}$) if and only if $w$ is beyond $v'$ from the
point of view of $v$.  It follows that the admissible monomials at
$v'$ all have the same $\nu$--weight except for the near monomial for
$v$, which has higher $\nu$--weight. Hence, the ideal defining the
associated graded ring $R$ contains polynomials obtained from our
equations by setting coefficients of near monomials equal to zero. We
will show these generate an ideal whose quotient is an integral domain
of dimension 2, hence yield the full associated graded.

For convenience of notation we will take $e=e_n$ for this proof. We
may assume (see Section \ref{sec:equations}) that the equations
associated to the node $v$ are
$$
M_{e_i}+a_iM_{e_{n-1}}+b_iM_{e_n}=0, \quad i=1,\dots, n-2.
$$
For each $j=1,\dots,n$, let $\Delta_j$ be the tree cut off away from
$v$ at the midpoint of $e_j$ (so $\Delta_n=\Delta_e^R$). By the
Scholium to Theorem \ref{th:cisg}, the equations for the associated
graded $R$ that correspond to nodes in $\Delta_j$ give a complete
intersection description of the monomial curve
$\C[X_j^{\sg(\Delta_j)}]$.  Let $\phi_j\colon\C[z_w, w\text{ a leaf of
  }\Delta_j]\to \C[X_j^{\sg(\Delta_j)}]$ be the corresponding
homomorphism. Then $\phi_j(M_{e_j})=c_jX^{p_j}$ for some $c_j\in\C^*$.
Together these homomorphisms $\phi_j$ give a homomorphism $\phi$ of
$R$ to the \BCI{} defined by the equations
\begin{equation}\label{eq:bci}
c_iX_i^{p_i}+c_{n-1}a_iX_{n-1}^{p_{n-1}}+c_nb_iX_n^{p_n}=0,\quad
i=1,\dots,n-2\,.
\end{equation}
Let $z_1,\dots,z_k$ be the variables corresponding to nodes of
$\Delta$ in $\Delta_e^L$ and let\break $z_{k+1},\dots,z_N$ be the remaining
variables, corresponding to nodes in $\Delta_e^R$. The graded
equations corresponding to nodes in $\Delta_e^L$ are equations for the
complete intersection curve defined by $\Delta_e^L$ except for
additional terms $b_iM_{e_n}$ (in the equations corresponding to node
$v$).  The procedure of subsection \ref{subsec:normalform} to put a
monomial in normal form will therefore change a monomial $M$ in the
variables $z_1,\dots,z_k$ into a linear combination of monomials of
the form $M'M_{e_n}^\alpha$, $\alpha\ge0$, with $M'\in\mathcal M_e^L$.
Thus, given any monomial in $z_1,\dots, z_N$, we first apply the
graded equations corresponding to nodes in $\Delta_e^L$ to put
anything involving $z_1,\dots,z_k$ in $\mathcal M_e^L$--normal form (at
the expense of adding factors $M_{e_n}$), and then apply the graded
equations corresponding to nodes in $\Delta_e^R$ to put anything
involving $z_{k+1},\dots,z_N$ into $\mathcal M_e^R$--normal form. It
follows that the set $\mathcal M_e^L \mathcal M_e^R$ is a $\C$--spanning
set for the graded ring $R$.  On the other hand, one can
check that the set
$$\{\phi(M_1M_2):M_1\in \mathcal M_e^L, M_2\in \mathcal M_e^R\}\subset
\C[X_1,\dots,X_n]/(\text{relations (\ref{eq:bci})})$$
is linearly
independent (we will not give a detailed proof of this, since it is
immediate in the case below to which we apply this proposition). Hence
$\mathcal M_e^L \mathcal M_e^R$ is a monomial basis for $R$. Moreover,
since $\phi$ is birational, $\phi$ is the normalization of $R$.
Finally, the calculation of $k_{v}$ is straightforward, using either
Proposition \ref{prop:k} or Proposition \ref{prop:K}.
\end{proof}

Note that the monomial basis given by the above proposition depends on
the choice of edge and also on the ordering of the variables. Although
the proposition gives the same monomial basis for the valuations
corresponding to the two ends of the edge, if we take a different node
we will have to take a different edge and will in general get a
different monomial basis.  However, to apply this proposition to
compute the geometric genus of $(X,p)$ we shall need the same monomial
basis for all the valuations. This turns out to be possible for
splice diagrams of the type
$$\def\lab#1{#1}\def\labb#1{#1} \splicediag{12}{16}{
  \Circ&& && && & &&  && &&\Circ
  \\ \hbox to 0pt{\hss$\Delta~=$}& \Vdots&\Circ\lineto[ull]_(.5){\lab{p_1}}
  \lineto[dll]^(.5){\lab{p_{k_1}}}
  \lineto[rr]^(.35){\labb{q_1}}^(.65){\labb{q'_1}}&&
  \Circ\lineto[ddl]_(.5){\lab{p_{k_1+1}}}
  \lineto[ddr]^(.5){\lab{p_{k_2}}}
  \lineto[r]^(.35){\labb{q_2}}&\dashto[r]&~ &
  \dashto[r]&\lineto[r]^(.35){\labb{q'_{r-1}}}&
  \Circ\lineto[ddl]_(.5){\lab{p_{k_{r-1}+1}}}
  \lineto[ddr]^(.5){\lab{p_{k_{r}}}}
  \lineto[rr]^(.35){\labb{q_{r}}}^(.65){\labb{q'_r}}&&
  \Circ\lineto[urr]^(.5){\lab{p_{N}}}
  \lineto[drr]_(.5){\lab{p_{k_r+1}}}
  &\Vdots\\
  \Circ&&&&\dots&& & && \dots&&&&\Circ\\
  &&&\Circ&&\Circ&&&\Circ&&\Circ}
$$
We number the nodes and edges of this diagram $v_0,\dots,v_r$ and
$e_1,\dots,e_r$ from left to right.  The valuation for node $v_i$ will
be denoted $\nu_i$.

For the edge $e=e_i$ joining nodes $v_{i-1}$ and $v_{i}$ we divide the
variables $z_1,\dots,z_{N}$ into two groups, ordered as follows:
\begin{gather*}
z_{k_{i}}, z_{k_{i}-1},\dots, z_1, \\
z_{k_{i}+1}, z_{k_i+2},\dots, z_{N}
\end{gather*}
We apply the above proposition for this particular edge $e$.
Example \ref{ex:monomialbasis} gives the monomial bases
\begin{align*}
\mathcal M_e^L&=\{z_1^{\alpha_1}\dots
z_{k_i}^{\alpha_{k_i}}:
0\le\alpha_i<p_i\text{ for }i=2,\dots,k_i\}\\
\mathcal M_e^R&=\{z_{k_i+1}^{\alpha_{k_i+1}}\dots
z_{N}^{\alpha_{N}}:
0\le\alpha_i<p_i\text{ for }i=k_i+1,\dots,N-1\}
\end{align*} for the two semigroups in question, so we get:

\begin{lemma}\label{lem:7.5}
  For each valuation $\nu_i$ of the above $\Delta$,
  $$\mathcal M:=\{z_1^{\alpha_1}\dots
  z_N^{\alpha_N}: 0\le\alpha_i<p_i\text{ for }i=2,\dots,N-1\}$$
  is a
  monomial basis for the associated graded ring $R$.\qed
\end{lemma}

We continue to consider the edge $e=e_i$ of $\Delta$ with left end
node $v=v_{i-1}$.  We can consider $\Delta^L_e$ and $\Delta^R_e$ also
as splice diagrams, and then $\Delta$ is the result of splicing them
at their root leaves.
\begin{theorem}\label{th:2n}
  The geometric genus $p_g(\Delta)$ of the splice type singularity
  determined by $\Delta$ is given inductively by
$$(1/4)C^LC^R+ p_{g}(\Delta_e^L) + p_g(\Delta_e^R)  $$
where $C^L$, $C^R$ are the conductors of the semigroups
$\sg(\Delta_e^L,w')$ and $\sg(\Delta_e^R,w')$
\end{theorem}
The following is a corollary of this and of Theorem \ref{th:equiv} in
the next section.
\begin{corollary}\label{cor:7.7}
  The Casson Invariant Conjecture holds for the splice type
  singularity determined by the above splice diagram $\Delta$
\end{corollary}
\begin{proof}[Proof of Corollary]
  The above theorem reduces this to an induction. The induction step
  is provided by Theorem \ref{th:equiv}, since $C^L$ and $C^R$ are the
  Milnor numbers of the knots corresponding to the root leaves of
  $\Delta_e^L$ and $\Delta_e^R$ (Theorem \ref{th:milnor}).
\end{proof}
\begin{proof}[Proof of Theorem \ref{th:2n}]
The canonical ideal of the singularity $(X,p)$ consists of those $f$
for which $\nu_j(f)\geq k_{v_j}$ for $j=0,\dots,r$.  Using the
linearly independent monomials of the above lemma, the geometric genus
thus equals the number of elements $M$ of
  $$\mathcal M=\{z_1^{\alpha_1}\dots
    z_N^{\alpha_N}: 0\le\alpha_i<p_i\text{ for }i=2,\dots,N-1\}$$
satisfying
\begin{equation}
  \label{eq:4}
\nu_i(M)< k_{v_i}\qquad\text{for some $i=0,\dots,r$.}
\end{equation}
We will call this condition ``\emph{condition $K(v_i)$}.'' So we want
to count the $M\in\mathcal M$ for which condition $K(v)$ holds for some
node $v$.

Let $e=e_i$.  For a monomial $M=z_1^{\alpha_1}\dots z_N^{\alpha_N}$,
write $M=M_LM_R$ with $M_L=z_1^{\alpha_1}\dots m_{k_i}^{\alpha_{k_i}}$
and $M_R=z_{k_i+1}^{\alpha_{k_i+1}}\dots z_N^{\alpha_N}$.  The
monomial $M$ is in $\mathcal M$ if and only if $M_L$ and $M_R$ are
normal form monomials for the semigroups $\sg(\Delta_e^L)$ and
$\sg(\Delta_e^R)$.

Denote the nodes at the left and right end of $e=e_i$ by
$v=v_{i-1}$ and $v'=v_i$ and the associated valuations by
$\nu=\nu_{i-1}$ and $\nu'=\nu_i$.  Denote
$$
  \ell_e(M_L) :=\sum_{j=1}^{k_i}\alpha_j\ell_{w'w_j}\,,\qquad
  \ell_e(M_R) :=\sum_{j=k_i+1}^N \alpha_j\ell_{w'w_j}
  $$
  where $w'$ is the root vertex of $\Delta_e^L$ or $\Delta_e^R$
  and $\ell_{w'w_j}$ is computed in $\Delta_e^L$ or $\Delta_e^R$.
  (Thus $\ell_e(M_L)$ and $\ell_e(M_R)$ are the values in the
  semigroups $\sg(\Delta_e^L)$ and $\sg(\Delta_e^R)$ corresponding to
  the monomials $M_L$ and $M_R$.) Then
$$\nu(M)=d_{ve}\ell_e(M_L)
+({d_{v}}/{d_{ve}})\ell_e(M_R).
$$
By Proposition \ref{prop:assocgraded} condition $K(v_{i-1})$
can thus be written
\begin{equation}
  \label{eq:7}
  d_{ve}(\ell_e(M_L)-C^L+1) +({d_{v}}/{d_{ve}})(\ell_e(M_R)-C^R+1)<0
\end{equation}
By symmetry, condition $K(v_i)$ can be written
\begin{equation}
  \label{eq:8}
({d_{v'}}/{d_{v'e}})(\ell_e(M_L)-C^L+1)
+d_{v'e}(\ell_e(M_R)-C^R+1)<0.
\end{equation}
Denote
$$  X_{i}:=\ell_e(M_L)-C^L+1\qquad Y_{i}:=\ell_e(M_R)-C^R+1\,,$$
so (\ref{eq:7}) and (\ref{eq:8}) can be written
\begin{equation}
  \label{eq:6}
  \begin{array}{rlr}
    K(v_{i-1}):&& d_{ve}X_i +({d_{v}}/{d_{ve}})Y_i<0,\\[2mm]
    K(v_i):&&({d_{v'}}/{d_{v'e}})X_i+d_{v'e}Y_i<0.
  \end{array}
\end{equation}
Note that $X_i\ne0$ since
$\ell_e(M_L)$ is in the semigroup $\sg(\Delta_e^L)$ with conductor
$C^L$. Similarly $Y_i\ne0$.  We will count the monomials $M\in
\mathcal M$ that satisfy condition $K(v)$ for some node $v$ by
subdividing into the following cases.
\begin{enumerate}
\item $X_i<0$ and $Y_i<0$ (so $K(v_{i-1})$ and $K(v_i)$ hold),
\item $Y_i>0$ and $K(v_j)$ holds for some $j\le i-1$,
\item $X_i>0$ and $K(v_j)$ holds for some $j\ge i$,
\item $Y_i>0$ and $K(v_j)$ holds for some $j\ge i$ and fails for all
  $j\le i-1$.
\item $X_i>0$ and $K(v_j)$ holds for some $j\le i-1$ and fails for all
  $j\ge i$,
\end{enumerate}
These cases cover all possibilities. We shall show that cases (4) and
(5) are empty and that Cases (1), (2), (3) are mutually exclusive and
lead to the three terms on the right in the theorem.

(1)\qua The number of monomials $M_L$ in normal form with $\ell_e(M_L)<
C^L-1$ is the number of elements bounded by $C^L$ in the semigroup
$\sg(\Delta_e^L)$. This is exactly $C^L/2$. Similarly for $M_R$, so
the set of $M\in\mathcal M$ with both $\ell_e(M_L)-C^L+1<0$ and
$\ell_e(M_R)-C^R+1<0$ contributes the $(1/4)C^LC^R$ of the
theorem.

(2)\qua  The inequality $\ell_e(M_R)-C^R+1>0$ says $\ell_e(M_R)\ge C^R$, so
there exists a unique monomial $M_R$ in normal form with such a value
of $\ell_e(M_R)$.  That is, if we put $\alpha=\ell_e(M_R)-C^R$ then
there is no constraint on $\alpha\ge0$ for a corresponding $M_R$ to
exist.  Consider the monomials $M_LM_R$ and $M_Lz^\alpha$, which are
normal form monomials for the splice diagrams $\Delta$ and
$\Delta_e^L$ respectively.  A simple calculation, which we omit, shows
that $M_LM_R$ satisfies condition $K(v_j)$ for $\Delta$ with $j\le
i-1$ if and only if $M_Lz^\alpha$ satisfies $K(v_j)$ for $\Delta_e^L$.
Thus the monomials $M=M_LM_R$ satisfying (2) are in one-one
correspondence with the monomials that count $p_g(\Delta_e^L)$.

(3)\qua  By symmetry, these monomials count $p_g(\Delta_e^R)$.

(4)\qua  One calculates that
$$Y_i=(d_{v_i}/q'_iq_{i+1})Y_{i+1}+(d_{v_i}/q_i)
\sum_{j=k_i+1}^{k_{i+1}}\frac1{p_j}(\alpha_j+1-p_j)$$
(we are using the explicit weights $d_{v_ie}=q'_i$ etc.\ from the
picture of $\Delta$). Since $\alpha_j<p_j$ for
$j=k_i+1,\dots,k_{i+1}$, the sum on the right is non-positive so
$Y_i>0$ implies $Y_{i+1}>0$.  Thus, if we are in case (4) we can, by
increasing $i$ if necessary, assume that $Y_i>0$ and $K(v_j)$ holds
for $j=i$ and fails for $j=i-1$.  By (\ref{eq:6}) we then have
$$d_{v}/d_{ve}^2>-X_i/Y_i,\qquad d_{v'e}^2/d_v'<-X_i/Y_i\,.$$
Thus
$d_{v}/d_{ve}^2>d_{v'e}^2/d_v'$, whence
$d_vd_{v'}/(d_{ve}d_{v'e})>d_{ve}d_{v'e}$, contradicting the edge
determinant condition. Thus case (4) cannot happen, and by symmetry
the same holds for case (5).

It remains to show that cases (2) and (3) are mutually exclusive (Case
(1) is clearly disjoint from (2) and (3)). But if both (2) and (3)
hold then $K(v_i)$ must fail. Since $K(v_j)$ holds for some $j>i$ and
$Y_i>0$, the same argument as in (4) leads to a contradiction.
\end{proof}

\section{Milnor fibers}\label{sec:milnor fiber}

Suppose $\Sigma$ is the link of an isolated singularity at $0$ of
a complete intersection surface $X=f^{-1}(0)$, where $f$ is a map
$f=(f_1,\dots,f_{n-2})\colon(\C^n,0)\to(\C^{n-2},0) $. The
\emph{Milnor fiber} is the manifold $F:=f^{-1}(\delta)\cap
B(\epsilon)$ where $B(\epsilon)$ is a sufficiently small ball
about $0$ and $\delta$ is a general point of $\C^{n-2}$ very close
to the origin.  It is a smooth simply-connected piece of complex
surface with boundary $\Sigma$; it has a symmetric intersection
pairing on the second homology group, whose rank $b_{2}(F)$ is
usually denoted by $\mu$.  The Casson Invariant Conjecture says
that when $\Sigma$ is a homology sphere, $\sign(F)$ should equal
$8\lambda(\Sigma)$, where $\lambda(\Sigma)$ is the Casson
invariant.

The Casson invariant of $\Sigma$ is not hard to compute, and the
hurdle in confirming this conjecture for any particular example is to
understand $F$ well enough to compute $\sign(F)$. This has been done
for \BCIs{}. Thus, the conjecture could be
verified in this case---a one-node splice diagram
(see \cite{neumann-wahl90}, which also proves a few other cases).

Now suppose the equations $f_i(z_1,\dots,z_n)=0$, $i=1,\dots,n-2$, are
of splice type as above, corresponding to a splice diagram
$\Delta$. Thus the curve $z_j=0$ cuts out in $\Sigma$ the knot $K_j$
corresponding to the $j$-th leaf of $\Delta$. The link $(\Sigma, K_j)$
is a fibered link whose fiber $G_j$ can also be seen as the Milnor
fiber of the singularity at $0$ of the complete intersection curve
$(f_1,\dots,f_{n-2},z_j)^{-1}(0)$. The topology of this fiber and its
embedding in $\Sigma$ can be described by gluing together Milnor
fibers of appropriate links in the splice components of $\Sigma$ (see
\cite{eisenbud-neumann}).

We shall describe a conjectural iterative description of $F$ in terms
of the Milnor fibers of simpler complete intersection surface
singularities and fibers $G_j$ as above lying in their boundaries.

Thus consider $\Sigma$ as the splice
$\Sigma=\Sigma_1~\raise5pt\hbox{$\underline{K_1\quad K_2}$}~\Sigma_2$
of two homology spheres determined by cutting $\Delta$ at an edge to
form two rooted diagrams. It is easy to see that these two diagrams
$\Delta_1$ and $\Delta_2$ also satisfy the semigroup condition so
$\Sigma_1$ and $\Sigma_2$ are both complete intersection singularity
links given by equations of splice type.  They thus have
Milnor fibers, which we shall call $F_1$ and $F_2$, with $\partial
F_i=\Sigma_i$.

\def\Fconj{\overline F}
Let $G_1\subset \Sigma_1$ be the fiber for the knot $(\Sigma_1, K_1)$.
This is the Milnor fiber described in the paragraph before Theorem
\ref{th:milnor}, so it is topologically determined by the rooted
diagram $\Delta_1$ and $b_1(G_1)$ is computed as in that theorem.

We may push the embedding $G_1\to F_1$ inside $F_1$ by a normal
vector-field to obtain a proper embedding $G_1\to F_1$ (that is, an
embedding with $\partial G_1=G_1\cap\partial F_1$, transverse
intersection) and then extend to an embedding $G_1\times D^2\to F_1$
of a tubular neighborhood of $G_1$. We similarly construct an
embedding $D^2\times G_2\to F_2$.

Denote
$$F_1^o:=F_1-(G_1\times\interior D^2),\quad F_2^o:=F_2-(\interior
D^2\times G_i),\quad$$
so $\partial F_1^o$ is the union of $G_1\times
S^1$ and the exterior (complement of an open tubular neighborhood) of
the knot $K_1\subset \Sigma_1$, and similarly for $\partial F_2^o$.

\begin{conjecture}[Milnor Fiber Conjecture]\label{conj:milnor fiber}
  $F$ is homeomorphic to the result $\Fconj$ of pasting:
  $$\Fconj:=F_1^o\cup_{G_1\times S^1} (G_1\times G_2)\cup_{S^1\times
    G_2}F_2^o,$$
  where we identify $G_1\times S^1$ with
  $G_1\times\partial G_2$ and $S^1\times G_2$ with $\partial G_1\times
  G_2$.
\end{conjecture}

By Milnor \cite{milnor-book} and Hamm \cite{hamm}, $F$, $F_1$, $F_2$
are simply connected $4$--manifolds which are homotopy equivalent to
$2$--complexes and thus have reduced homology only in dimension $2$.
We show that $\overline F$ has the nice properties we would like $F$
to have.
\begin{theorem} \label{th:homol}   $\partial\Fconj=\Sigma$ and
$\Fconj$ is  simply connected and
  homotopy equivalent to a $2$--complex. Moreover,
\begin{align*}
  H_2(\Fconj)&\isom H_2(G_1\times G_2)\oplus H_2(F_1)\oplus
H_2(F_2)\\&=\bigl(H_1(G_1)\otimes H_1(G_2)\bigr)\oplus H_2(F_1)\oplus
H_2(F_2).
\end{align*} with maps induced by inclusions, so
  $$\sign(\Fconj)=\sign(F_1)+\sign(F_2)\,.$$
\end{theorem}
\begin{corollary}\label{cor:milnor fiber}
  The Milnor Fiber Conjecture (Conjecture \ref{conj:milnor fiber})
  implies the Casson Invariant Conjecture for complete intersection
  singularities of splice type.
\end{corollary}
\begin{proof}
  The theorem and Conjecture \ref{conj:milnor fiber} imply that
  signature of Milnor fiber is additive under splicing. The Casson
  invariant is additive for splicing.  The Casson Invariant Conjecture
  is known for \BCIs{}
%equations of \Brieskorn{} type
(the one-node case).
\end{proof}
\begin{proof}[Proof of Theorem \ref{th:homol}]
  The fact that $\partial\Fconj=\Sigma$ is immediate from the
  construction.  For the rest of this proof it is convenient to have a
  different description of $\Fconj$.

Consider $G_i$ embedded in $\Sigma_i$ and let $N_i\subset
\Sigma_i=\partial F_i$ be a tubular neighborhood of $G_i$ in
$\Sigma_i$, so $N_i\isom G_i\times I$.  Note that $\partial (G_1\times
G_2)=(G_1\times K_2)\cup(K_1\times G_2)$, so we can also embed $N_1$
in $\partial (G_1\times G_2)$ as $G_1\times I\subset G_1\times K_2$,
and similarly for $N_2$. We claim:
\begin{equation}
  \label{eq:5}
\Fconj\cong F_1\cup_{N_1}(G_1\times G_2)\cup_{N_2}F_2.
\end{equation}
Indeed, to turn our previous description of $\Fconj$ into this one,
connect the
proper embedding $G_i\subset F_i$ to the
embedding $G_i\subset \partial F_i$
%to the previous embedding $G_i\times D^2\subset F_i$
by a ``strip'' $G_i\times I$ and remove a
tubular neighborhood of this strip from $F_i^o$ and glue it onto
$G_1\times G_2$ instead. The result of removing it from $F_i^o$ is
something homeomorphic to $F_i$, while, when glued to $G_1\times
G_2$ it is just a collar on part of the boundary and does not change
the homeomorphism type of $G_1\times G_2$.

Consider, therefore, $\Fconj$ as in equation (\ref{eq:5}).
By shrinking slightly the regions $N_i$ along which the $F_i$ are
glued to $G_1\times G_2$ we can make them disjoint in
$\partial(G_1\times G_2)$ without changing the homotopy type (or even
homeomorphism type) of $F_1\cup_{N_1}(G_1\times G_2)\cup_{N_2}F_2$.
Then $(G_1\times G_2)\cap (F_1\cup F_2)$ consists of the disjoint
union of $N_1$ and $N_2$.  The Meyer-Vietoris sequence for the
decomposition $(G_1\times G_2)\bigcup(F_1\cup F_2)$ then easily yields
that the inclusions induce an isomorphism
$$  H_2(\Fconj)\isom H_2(G_1\times G_2)\oplus H_2(F_1)\oplus
H_2(F_2)
$$ as desired.

The fact that $\Fconj$ is simply connected is an easy application of
the Van Kampen theorem. The fact that $\Fconj$ is homotopy equivalent
to a $2$--complex can be seen by replacing $G_1$ and $G_2$ by
one-dimensional spines $S_1$ and $S_2$ say, replacing $F_1$ and $F_2$ by
2-dimensional spines $T_1$ and $T_2$, and then gluing $S_1\times S_2$
to $T_1$ and $T_2$ by means of mapping cylinders of appropriate
maps $S_i\to T_i$.
\end{proof}

Recall that the geometric genus $p_g(X,o)$
of a singularity is $\dim H^1(Y,\mathcal O)$, where $Y \to X$
denotes a resolution of the singularity.  In general, it is not
topologically determined by the link of $X$, but the \CIC{}
%(second version)
says that it should be for
%Gorenstein
complete intersection singularities with
homology sphere links. The following theorem says what the \CIC{}
implies about the behavior of various invariants under splicing. Item
(\ref{it:betti}) of this theorem provided part of the motivation for
the above construction of $\Fconj$ for the Milnor Fiber Conjecture.

\begin{theorem}\label{th:equiv}
  Let $X$ be a complete intersection with homology sphere link, with
  Milnor fiber $F$; and suppose its link is spliced from links of two
  singularities $X_{1},X_{2}$, with Milnor fibers $F_{1},F_{2}$.
  Assume the Casson Invariant Conjecture for $X_{1}$ and $X_{2}$.
  Then the following statements are equivalent:
  \begin{enumerate}
  \item The Casson Invariant Conjecture holds for $X$.
  \item We have $\sign(F)=\sign(F_1)+\sign(F_2)$.
  \item With $G_1$, $G_2$ as above, we have
    $b_2(F)=b_2(F_1)+b_2(F_2)+b_1(G_1)b_1(G_2)\,,$ where $b_i$ is
    Betti number.\label{it:betti}
  \item The geometric genus satisfies $p_g(X)=p_g(X_{1})+p_g
    (X{_2})+\frac14b_1(G_1)b_1(G_2)$.
  \end{enumerate}
  Moreover, these invariants of $X$ are then topologically determined
  by the link.
\end{theorem}
\begin{proof}
The equivalence of (1) and (2) has already been discussed, so we prove
the equivalence of (2), (3), and (4).

Formulas of H. Laufer and A. Durfee imply that the geometric genus of
X, and the signature and second Betti number $\mu$ of the Milnor
fiber, are explicitly related by topological invariants of the link
(see, eg, \cite{wahl}.)  Let $Y\rightarrow X$ denote a good
resolution, and $c_1^2$ and $c_2$ the
characteristic Chern numbers of $Y$ (also known as $K\cdot
K$ and $\chi(Y)$, where $\chi$ is topological Euler
characteristic).  Then these Chern numbers are
determined by the resolution dual graph, and their sum $c_{1}^{2}+c_{2}$ is
independent of the resolution, hence depends only on the link.
We define $$C(\Delta)=c_{1}^{2}+c_{2}-1,$$ the notation indicating
that this number depends only on the splice diagram $\Delta$.  Then
the aforementioned formulas may be written
\begin{align*}
\mu&=12p_g+C(\Delta)&\text{(Laufer)}\\
3\sign(F)&=-2\mu-C(\Delta)&\text{(Durfee)}
\end{align*}
(In general Durfee's formula has an extra $3b_1(Y)$ on the right,
which vanishes in our case.)  Eliminating $\mu$, these formulas imply
$$\sign(F)=-8p_g -C(\Delta)\,,$$
proving the equivalence of the two
formulations of the \CIC{} in the Introduction (for complete
intersections).
Moreover, it follows that the equivalence of (2), (3), and (4) of
Theorem \ref{th:equiv} reduce to the formula of the following theorem,
which will therefore complete the proof.
\end{proof}
\begin{theorem}\label{th:Cdelta}
In the above notation, even if $\Delta$ does not satisfy the semigroup
condition we have
$$C(\Delta)-C(\Delta_1)-C(\Delta_2)=-2b_1(G_1)b_1(G_2).$$
\end{theorem}
This theorem involves computing $c_1^2$ and $c_2$ of the resolution in
terms of the splice diagram, which is of interest in its own right, so
we devote the next section (Section \ref{sec:canonical}) to its proof.

If any one of the analytic invariants $\sign(F)$, $\mu$, and $p_g(X)$
is a topological invariant, then they all are, by the above formulas.
The Casson Invariant Conjecture gives a topological description of
$p_g$ and $\sign(F)$.

Suppose $(X,o)$ is a
%Gorenstein
complete intersection surface singularity whose homology sphere link
has one node; thus, its link is $\Sigma(p_{1},\cdots,p_{n})$.  The
Casson Invariant Conjecture for $X$ is equivalent to the assertion
that $p_{g}(X)=p_{g}(V(p_{1},\cdots,p_{n}))$.  But this latter
condition is well-known to be equivalent to the statement that $X$
admits an equisingular, simultaneous resolution degeneration to $V$
(see, eg, \cite{tomari-watanabe} (6.3) for a convenient proof).  In
other words, we could conclude that $X$ is a splice type singularity,
as was mentioned in the Introduction.  We suspect a similar result is
true in the general case.  But, even in case $X$ is a hypersurface
singularity with link $\Sigma(p,q,r)$, we do not know a proof.  As we
mentioned in the Introduction, there are a few very non-trivial cases
worked out by A N\'emethi \cite{nemethi}.

\section{Canonical divisor of a resolution}\label{sec:canonical}

This section is devoted to proving Theorem \ref{th:Cdelta}. We start by
computing the rational canonical divisor for an arbitrary resolution
of an isolated surface singularity.

Suppose we have a good  resolution of an isolated surface singularity.
Denote the exceptional curves by  $E_i$, $i=1,\ldots,n$.  For each $i$
let $\delta_i$  be  the number  of intersection  points of $E_i$  with
other $E_j$'s and  let $E_i^0$  be $E_i$  with these  intersection  points
removed.    Denote  $\chi_i=\chi(E_i^0)=\chi(E_i)-\delta_i$ ($\chi$ is
Euler characteristic).

Let $K$ be the (rational) canonical divisor, defined by the adjunction
formula
$$K\cdot E_i= -\chi(E_i)-E_i\cdot E_i.$$
$$D := -K-E,\quad\hbox{where }E=\sum_{i=1}^n E_i\leqno{\rm Let}$$
$$D=\sum_{i=1}^n k_i E_i.\leqno{\rm and\ suppose}$$
Then the adjunction formula becomes
$D.E_j = \chi(E_j)-\delta_j = \chi_j$, so
$$k_i=-\sum \ell_{ij}\chi_j,\quad\hbox{ where } (\ell_{ij})=(-E_i\cdot
E_j)^{-1}\text{ (matrix inverse).}$$
$$K=-D-E=\sum_i(-k_i-1)E_i=\sum_i\Bigl(\sum_j
\ell_{ij}\chi_j-1\Bigr)E_i\leqno{\rm Now}$$
so
\begin{align*}
K\cdot K&=\biggl(\sum_i\Bigl(\sum_j \ell_{ij}\chi_j-1\Bigr)E_i\biggr)\cdot\biggl(\sum_k\Bigl(\sum_l \ell_{kl}\chi_l-1\Bigr)E_k\biggr)\\
&=\sum_{i,j,k,l}\ell_{ij}\ell_{kl}(E_i\cdot E_k)\chi_j\chi_l
 - \sum_{i,k,l}\ell_{kl}(E_i\cdot E_k)\chi_l
 - \sum_{i,j,k}\ell_{ij}(E_i\cdot E_k)\chi_j\\
&\qquad\qquad + \sum_{i,k}E_i\cdot E_k\\
&=-\sum_{j,l}\ell_{jl}\chi_j\chi_l
 + \sum_i\chi_i  + \sum_i\chi_i
 +\Bigl(\sum_i E_i\cdot E_i + 2 \sum_{i<j}E_i\cdot E_j\Bigr) \\
&= -\sum_{i,j}\ell_{ij}\chi_i\chi_j +2 \chi(\bigcup_i E_i) +\sum_i E_i\cdot E_i\,.
\end{align*}
We note that our notation $k_i$ and $\ell_{ij}$ is consistent with the
notation in the Appendix (Section \ref{sec:splicing}).
Summarizing:%, using the notation of the proof of Theorem \ref{th:equiv}:
\begin{proposition} \label{prop:K}
For any good resolution of an isolated surface
  singularity   the divisor $D=-K-E$ is given by
$$D=\sum k_i E_i\quad\text{with}\quad k_i=-\sum_j \ell_{ij}\chi_j$$
where $(\ell_{ij})=(-E_i\cdot
E_j)^{-1}$ (matrix inverse).
Also,
$$c_1^2+c_2=-\sum_{i,j}\ell_{ij}\chi_i\chi_j + 3c_2 +\sum E_i\cdot
E_i\,.$$
\end{proposition}

To apply this to prove Theorem \ref{th:Cdelta} we now restrict to the
case of a singularity with homology sphere link given by a splice
diagram $\Delta$. Then $\chi_i=2-\delta_i$ which vanishes except at
nodes and ends of the plumbing graph, so we only need to know
$\ell_{ij}$ when $i$ and $j$ index nodes or ends.  By Theorem
\ref{th:props}, this is as follows. If $i\ne j$ then $\ell_{ij}$ is
the product of the splice diagram weights adjacent but not on the path
from $i$ to $j$ in $\Delta$.  If $i=j$ then:
\begin{itemize}
\item If $i$ is a node then $\ell_{ii}$ is the product of weights at that node.
\item If $i$ is a leaf adjacent to a node with weights
  $p_0,\ldots,p_n$ with $p_0$ on the edge to $i$ then $\ell_{ii}=
  \lceil p_1\ldots p_n/p_0\rceil$.% (integer part).
\end{itemize}

A simple matrix calculation shows that $-\ell_{ii}$ is the weight one
would have to put on a new vertex attached to vertex $i$ by a new
edge, to get an extended plumbing diagram of determinant $0$ (this is
computed in \cite{eisenbud-neumann} and gives an alternative proof of
the description of $\ell_{ii}$).

Now let $C(\resgraph)$ denote $c_1^2+c_2-1$ computed for a plumbing graph
$\resgraph$.  We want to compute the effect of splicing on $C$.  So
suppose that the splice diagram $\Delta$ is the result of splicing
diagrams $\Delta_1$ and $\Delta_2$ and let $\resgraph$, $\resgraph_1$,
and $\resgraph_2$ be the resolution graphs for these splice diagrams.
Let $I_1$ and $I_2$ be index sets for the nodes and
leaves of $\resgraph_1$ and $\resgraph_2$ with $0\in I_1$ and $1\in I_2$
representing the leaves at which we splice.
%Then $I=I_1\bigcup I_2 -
%\{0,1\}$ is the index set for nodes and leaves of $\resgraph$.

Let $\overline\resgraph_1$ be the result of extending $\resgraph_1$ at
vertex $0$ by a vertex with weight $-\ell_{00}(\resgraph_1)$ and similarly
for $\overline\resgraph_2$.  Let $\overline\resgraph$ be the result of
attaching $\overline\resgraph_1$ to $\overline\resgraph_2$ by an edge
joining the new vertices.  Then in \cite{eisenbud-neumann} it is shown
that $\resgraph$ results from $\overline\resgraph$ by a sequence of
$(-1)$--blow-downs of vertices of valency 2 followed by one
$0$--absorption.  Suppose the number of blow-downs is $r$.  Then the
blow-downs and $0$--absorption remove (r+2) vertices, so
$$c_2(\resgraph)=c_2(\overline\resgraph)-(r+2)=
c_2(\resgraph_1)+c_2(\resgraph_2)-r-1.$$
Moreover each blow-down increases $\sum E_i\cdot E_i$ by $3$ and the
$0$--absorption does not change it, so
$$\sum_{\resgraph}E_i\cdot E_i=\sum_{\resgraph_1}E_i\cdot E_i - \ell_{00}
+\sum_{\resgraph_2}
E_i\cdot E_i - \ell_{11}\quad +\,\,3r\,.$$
Thus
\def\backspace#1pt{\hbox to -#1pt{}}
\begin{align*}
C(\resgraph)
&= -\sum_{i,j\in I}\ell_{ij}\chi_i\chi_j\quad +\quad 3c_2({\resgraph})
\quad +\quad \sum_{\resgraph}E_i\cdot E_i\quad-\quad 1\\
&= -\Bigl(\sum_{i,j\in I_1}\ell_{ij}\chi_i\chi_j - \ell_{00} - 2\backspace10pt
\sum_{i\in I_1-\{0\}}\backspace7pt \ell_{0i}\chi_i
+ \sum_{i,j\in I_2}\ell_{ij}\chi_i\chi_j - \ell_{11} - \\
&\qquad\qquad - 2\backspace10pt \sum_{j\in I_2-\{1\}}
\backspace7pt \ell_{1j}\chi_j
+ 2\backspace10pt \sum_{i\in I_1-\{0\}\atop
j\in I_2-\{1\}}\backspace7pt \ell_{ij}\chi_i\chi_j\Bigr) \\
&\quad + 3\bigl(c_2({\resgraph_1}) + c_2({\resgraph_2}) -r-1\bigr)\\
&\quad +\sum_{\resgraph_1}E_i\cdot E_i +\sum_{\resgraph_2}E_i\cdot E_i
-\ell_{00}- \ell_{11}+3r\quad-\quad1\\
&=
-\sum_{i,j\in I_1}\ell_{ij}\chi_i\chi_j + 2\backspace10pt \sum_{i\in
  I_1-\{0\}}\backspace7pt \ell_{0i}\chi_i
- \sum_{i,j\in I_2}\ell_{ij}\chi_i\chi_j \qquad\qquad\\
&\qquad\qquad + 2\backspace10pt \sum_{j\in I_2-\{1\}}
\backspace10pt \ell_{1j}\chi_j
- 2 \backspace7pt \sum_{i\in I_1-\{0\}}\backspace7pt \ell_{i0}\chi_i
\backspace7pt \sum_{j\in I_2-\{1\}}\backspace7pt \ell_{j_1}\chi_j \\
&\quad + 3c_2(\resgraph_1) + 3c_2(\resgraph_2) \\
&\quad +\sum_{\resgraph_1}E_i\cdot E_i +\sum_{\resgraph_2}E_i\cdot E_i\quad-\,\,4\,,
\end{align*}
where the last equality uses the fact that for $i\in I_1-\{0\}$ and
$j\in I_2-\{1\}$ one has $\ell_{ij}=\ell_{i0}\ell_{1j}$.

The above simplifies to
\begin{align*} C(\resgraph)&=C(\resgraph_1)+C(\resgraph_2)-2\Bigl(-\backspace10pt\sum_{i\in
    I_1-\{0\}}\backspace10pt \ell_{0i}\chi_i\,+1\Bigr)
  \Bigl(-\backspace10pt\sum_{j\in
    I_2-\{1\}}\backspace10pt \ell_{1j}\chi_j\,+1\Bigr)\\
  &=C(\resgraph_1)+C(\resgraph_2)-2\mu(\resgraph_1,0)\mu(\resgraph_2,1) ,
\end{align*}
where $\mu(\resgraph_i,i)$ is the Milnor number for the knot represented
by vertex $i$ in the homology sphere represented by $\resgraph_i$, that
is, the first Betti number of its fiber (it is a basic result of
\cite{eisenbud-neumann} that $\sum_{i\in I_1-\{0\}}(\ell_{0i}\chi_i)$ is
the Euler characteristic of the fiber in question).  This completes
the proof of Theorem \ref{th:Cdelta}.\qed

\section{Plane curves and their cyclic covers}\label{sec:plane curves}

Let $(X,o)$ be a hypersurface singularity at the origin given by an
equation in the form $z^n+g(x,y)=0$ and suppose that its link is a
homology sphere. The Casson Invariant Conjecture was proved in this
case in \cite{neumann-wahl90} by a somewhat subtle calculation. In
this section we will show that the Milnor Fiber Conjecture (Conjecture
\ref{conj:milnor fiber}) holds for these singularities, giving a more
conceptual proof of the Casson Invariant Conjecture in this case. We
must first explain how these hypersurface singularities fit the format
of equations of splice type.
%Since these are instructive examples for our conjectures, we discuss
%some other aspects of them too.
\comment {Saveliev says he has a short proof, using the branched
  cover} In \cite{neumann-wahl90} we point out that if the link of
$z^n+g(x,y)=0$ is a homology sphere, then $g(x,y)=0$ defines an
irreducible plane curve singularity at the origin $o\in \C^2$. We
therefore need to start by discussing how plane curve singularities in
general, and irreducible plane curve singularities in particular, fit
into the framework of our conjectures.

\subsection{Non-minimal splice diagrams and plane curve singularities}
Theorem \ref{th:main1} gives a general sufficient condition for a knot
in a homology sphere to be realizable as the link of a germ
$(Y,o)\subset (X,o)$ of a curve cut out by a single equation in a
complete intersection surface.  This has content also for non-minimal
splice diagrams. For example, the splice diagram
$$
\splicediag{6}{30}{\\
&\Circ&&&\Circ\\
\Delta\quad=&&\Circ\lineto[ul]_(.25){2}\lineto[dl]^(.25){3}
&\Circ\lineto[dr]_(.25){2}\lineto[ur]^(.25){1}
\lineto[l]_(.2){3}_(.8){5}\\
&\Circ&&&\Circ
}$$
is a non-minimal version of
$$\splicediag{6}{30}{\\
&\Circ\\
&&\Circ\lineto[ul]_(.25){2}\lineto[dl]^(.25){3}\lineto[r]^(.25){5}
&\Circ\\
&\Circ
}
$$
so it represents the Seifert fibered homology sphere
$\Sigma(2,3,5)$ (Poincar\'e's dodecahedral space). The upper right
vertex of $\Delta$ represents a particular knot in this homology
sphere (a $(3,2)$--cable on the degree 5 fiber of $\Sigma(2,3,5)$).
Since $\Delta$ satisfies the semigroup condition, Theorem
\ref{th:main1} tells us that this knot in $\Sigma(2,3,5)$ is the link
of a complex curve singularity $(Y,o)$ cut out by a single equation in
$(V(2,3,5),o)$. In fact, the splice type equations for $\Delta$ can be
chosen as $z_1^2+z_2^3+z_3^5=0, z_1+z_3^2+z_4=0$, and the curve is
then cut out by $z_4=0$. Eliminating $z_4$, the curve is cut out by
the equation $z_1+z_3^2=0$ in
$V(2,3,5)=\{(z_1,z_2,z_3):z_1^2+z_2^3+z_3^5=0\}$.  \comment{also
  necessary?}

When $X$ is non-singular, that is, for a link of a plane curve
singularity, the next proposition implies that we can always do the
analogous thing.  That is, for any irreducible plane curve singularity
we will find strict splice type equations for $X$ ($=\C^2$) so that
the curve $Y$ cut out by a coordinate function has the topology of the
given plane curve. Corollary \ref{cor:znfxy} then says that the
original plane curve singularity is a higher weight deformation of the
one given by strict splice type equations.
\begin{proposition}
  The splice diagram of any plane curve singularity satisfies the
  semigroup condition.
\end{proposition}
\begin{proof}
  It is easy to see that the semigroup condition for the splice
  diagram of a reducible plane curve singularity follows from the
  semigroup condition for each of the subdiagrams for the irreducible
  branches of the plane curve.  Thus we may assume that the germ
  $(\C^2,Y,o)$ is an irreducible germ.  In this case the result is
  well known (see, eg, Teissier's appendix to \cite{zariski}) but we
  give a proof in our language for completeness. By
  \cite{eisenbud-neumann} the singularity is given by a splice diagram
  of the form:
$$%  \begin{equation}    \label{eq:sd0}
\splicediag{16}{24}{
    \Circ\lineto[r]^(.75){p_1}&\Circ\lineto[r]^(.25){1}^(.75){p_2}
    \lineto[d]^(.25){q_1}&\Circ\lineto[d]^(.25){q_2}\lineto[r]^(.25)1&
    \dotto[r]&\lineto[r]^(.75){p_k}&\Circ\lineto[d]^(.25){q_k}\ar[r]^(.25)1&\\
    &\Circ&\Circ&&&\Circ}
$$%  \end{equation}
  where $\gcd(p_i,q_i)=1$ for each $i$ and
  the positive edge determinant condition holds
  ($p_i>q_iq_{i-1}p_{i-1}$ for each $i>1$). Since this diagram may
  have arisen as a subdiagram of a diagram for a plane curve with
  several branches, we cannot assume that it is a reduced diagram, so
  some of the $q_j$ may equal $1$.

  The only non-trivial cases of the semigroup condition for this
  diagram are:
  $$p_{j+1}\in S_j:=\N\langle q_1q_2\dots q_{j},p_1q_2\dots
  q_{j},\dots,p_{j-1}q_{j}, p_{j}\rangle$$
  for each $j=1,\dots,k-1$.
  Since $p_{j+1}>p_jq_jq_{j+1}\ge p_jq_j$ it suffices to show that the
  conductor $\mu_j$ of this semigroup satisfies $\mu_j\le p_jq_j$.
  Proposition \ref{le:sg} of section \ref{sec:semigroup} implies
  $\mu_j=q_j(\mu_{j-1}-1)-p_j+p_jq_j+1$ (or use Theorem
  \ref{th:milnor} and its following paragraph). The desired inequality
  is now a trivial induction.
\end{proof}
This gives a new way to find an equation for a plane curve singularity
of given topology: start with the equations of splice type and then
eliminate variables to obtain an equation in $\C^2$.  To describe this
in detail, let us assign variables to the leaves of our splice diagram
as follows:
$$%\begin{equation}    \label{eq:sd1}
  \splicediag{6}{24}{z_0\,\,
    \Circ\lineto[r]^(.75){p_1}&\Circ\lineto[r]^(.25){1}^(.75){p_2}
    \lineto[ddd]^(.25){q_1}&\Circ\lineto[ddd]^(.25){q_2}\lineto[r]^(.25)1&
    \dotto[r]&\lineto[r]^(.75){p_k}&\Circ\lineto[ddd]^(.25){q_k}\lineto
   [r]^(.25)1&\Circ\,\,z_{k+1}\\ \\ \\
    &\Circ&\Circ&&&\Circ\\
  &z_1&z_2&&&z_k}
%\end{equation}
$$
The only admissible monomial for the outgoing edge to the right at the
$j$-th node is $z_{j+1}$.  Thus the general system of equations of
strict splice type can be written
\begin{align*}
  z_2&=a_{1}z_1^{q_1}+a_{0}z_0^{p_1}\\
  z_3&=a_{2}z_2^{q_2}+g_2(z_0,z_1)\\
  \dots&\quad\dots\quad \dots\\
  z_k&=a_{k-1}z_{k-1}^{q_{k-1}}+g_{k-1}(z_0,\dots,z_{k-2})\\
  z_{k+1}&=a_{k}z_{k}^{q_{k}}+g_{k}(z_0,\dots,z_{k-1}),
\end{align*}
where $g_j(z_0,\dots,z_{j-1})$ is a
multiple of an admissible monomial for the left edge at the $j$-th
node, that is, a monomial of
the form $z_0^{\alpha_0}\dots z_{j-1}^{\alpha_{j-1}}$ with
$$\alpha_0q_1\dots q_{j-1}+\alpha_1p_1q_2\dots
q_{j-1}+\dots+\alpha_{j-2}p_{j-2}q_{j-1}+\alpha_{j-1}p_{j-1}=p_{j}.$$

We now successively substitute each of the above equations into the next
to put them in the form:
\begin{align*}
  z_2&=a_{1}z_1^{q_1}+a_{0}z_0^{p_1}\\
  z_3&=a_2(a_1z_1^{q_1}+a_0z_0^{p_1})^{q_2}+g_2(z_0,z_1)=:f_2(z_0,z_1)\\
  \dots&\quad\dots\quad \dots\\
   z_{k+1}&=a_{k}f_{k-1}(z_0,z_1)^{q_{k}}+g_{k}(z_0,z_1,
  \dots,f_{k-2}(z_0,z_1))=:f_{k}(z_0,z_{1}).
\end{align*}
In terms of new coordinates, $x:=z_0$, $y:=z_1$, $
Z_2:=z_2-a_{1}z_1^{q_1}+a_{0}z_0^{p_1}$, $ \dots$, $
Z_{k}:=z_k-f_{k-1}(z_0,z_1)$, $ Z_{k+1}:=z_{k+1}-f_{k}(z_0,z_1)$ these
equations become
$$
Z_2=Z_3=\dots=Z_{k}=Z_{k+1}=0,
$$
so our surface is the $(x,y)$--plane. Our plane curve is the curve
cut out by the coordinate equation $z_{k+1}=0$ which is $f_{k}(x,y)=0$
in our new coordinates.  Thus, if we write $f=f_{k}$, the equation of
the plane curve is $f(x,y)=0$.

We now address what the Milnor Fiber Conjecture says for
this type of example. Our surface germ is a nonsingular point, and the
Milnor fiber for a non-singular point is a disk, so the conjecture
postulates a particular decomposition of $D^4$. Although it is rather
trivial, it will be needed in the discussion of hypersurfaces of the
form $z^n=g(x,y)$.  We will therefore reserve the notations $G_1$
etc.\ of Conjecture \ref{conj:milnor fiber} for that case and use
primes (as in $G_1'$ etc) to distinguish the ingredients involved in
the present discussion.

Suppose therefore that we have decomposed our splice diagram as the
splice of two diagrams:
$$\splicediag{16}{24}{
  \Circ\lineto[r]^(.75){p_1}&\Circ\lineto[r]^(.25){1}^(.75){p_2}
  \lineto[d]^(.25){q_1}&\Circ\lineto[d]^(.25){q_2}\lineto[r]^(.25)1&
  \dotto[r]&\lineto[r]^(.75){p_r}&\Circ\lineto[d]^(.25){q_r}\ar[r]^(.25)1&
  \quad&\Circ\ar[l]_(.25){p_{r+1}}\lineto[d]^(.25){q_{r+1}}\lineto[r]^(.25)1&
  \dotto[r]&\lineto[r]^(.75){p_k}&\Circ\lineto[r]^(.25)1\lineto[d]^(.25){q_k}&
  \Circ\\
  &\Circ&\Circ&&&\Circ&&\Circ&&&\Circ}
$$
The left diagram represents a plane curve whose Milnor fiber we
will denote by $G_1'\subset S^3=\partial D^4$. The right diagram is a
non-reduced diagram for the trivial knot in $S^3$ so its Milnor fiber
is $G_2'=D^2\subset S^3=\partial D^4$.

Let $(F_1')^o$ be the result of removing from $D^4$ a tubular
neighborhood of $G_1'$ pushed inside to a proper embedding
$G_1'\subset D^4$. Let $(F_2')^o$ be the result of removing a tubular
neighborhood of a proper embedding $D^2\subset D^4$. Note that
$(F_2')^o\cong S^1\times D^3$. The Milnor Fiber Conjecture says that
the result of the pasting:
\begin{equation}
  \label{eq:3}
%D^4\cup_{N_1'}G_1'\times D^2\cup_{N_2'}D^4
(F_1')^o\cup (G_1'\times D^2)\cup (F_2')^o
\end{equation}
should be $D^4$.  This is indeed clear, since, starting with
$(F_1')^o$, the first pasting clearly gives $D^4$ back, while the
second just pastes a collar onto a portion of the boundary of this
$D^4$.

\subsection{The hypersurface $z^n+g(x,y)=0$}\label{subsec:52}
As already mentioned, if the link of $z^n+g(x,y)=0$ is a homology
sphere, then $g(x,y)=0$ defines a plane curve singularity at $(0,0)\in
\C^2$ which is irreducible. Its splice diagram therefore has the form
$$\splicediag{16}{24}{
  \Circ\lineto[r]^(.75){p_1}&\Circ\lineto[r]^(.25){1}^(.75){p_2}
  \lineto[d]^(.25){q_1}&\Circ\lineto[d]^(.25){q_2}\lineto[r]^(.25)1&
  \dotto[r]&\lineto[r]^(.75){p_k}&\Circ\lineto[d]^(.25){q_k}\ar[r]^(.25)1&\\
  &\Circ&\Circ&&&\Circ}$$
where $\gcd(p_i,q_i)=1$ for each $i$ and
the positive edge determinant condition holds ($p_i>q_iq_{i-1}p_{i-1}$
for each $i>1$).  Moreover, given an irreducible plane curve
singularity as above, we showed in \cite{neumann-wahl90} that the
hypersurface singularity defined by $$z^n+g(x,y)=0$$
has homology
sphere link if and only if $n$ is relatively prime to all the $p_i$
and $q_i$, and the splice diagram for the link of this singularity is
then:
\begin{equation}
    \label{eq:sd}
\splicediag{16}{24}{
    \Circ\lineto[r]^(.75){p_1}&\Circ\lineto[r]^(.25){n}^(.75){p_2}
    \lineto[d]^(.25){q_1}&\Circ\lineto[d]^(.25){q_2}\lineto[r]^(.25)n&
    \dotto[r]&\lineto[r]^(.75){p_k}&\Circ\lineto[d]^(.25){q_k}\lineto
[r]^(.25)n&\Circ\\
    &\Circ&\Circ&&&\Circ}
\end{equation}
We now show that the splice diagram equations for this splice
diagram reduce to the equation
$z^n=f(x,y)$, with $f$ as in the previous subsection (Corollary
\ref{cor:znfxy} below shows that the original $z^n=g(x,y)$ is an
equisingular deformation of this).  We assign variables to the leaves
of the splice diagram (\ref{eq:sd}) as follows:
$$%  \begin{equation}
    \label{eq:sd2}
\splicediag{6}{24}{x=z_0\,\,
    \Circ\lineto[r]^(.75){p_1}&\Circ\lineto[r]^(.25){n}^(.75){p_2}
    \lineto[ddd]^(.25){q_1}&\Circ\lineto[ddd]^(.25){q_2}\lineto[r]^(.25)n&
    \dotto[r]&\lineto[r]^(.75){p_k}&\Circ\lineto[ddd]^(.25){q_k}\lineto
[r]^(.25)n&\Circ\,\,z\\ \\ \\
    &\Circ&\Circ&&&\Circ\\
&y=z_1&z_2&&&z_k}
%  \end{equation}
$$
  The only admissible monomial for the outgoing edge to the right at
  the $j$-th node is $z_{j+1}$ if $j<k$ and $z^n$ if $j=k$.
  Thus the general system of equations of strict splice type can be
  written
\begin{align*}
z_2&=a_{1}z_1^{q_1}+a_{0}z_0^{p_1}\\
z_3&=a_{2}z_2^{q_2}+g_2(z_0,z_1)\\
\dots&\quad\dots\quad \dots\\
z_k&=a_{k-1}z_{k-1}^{q_{k-1}}+g_{k-1}(z_0,\dots,z_{k-2})\\
z^n&=a_{k}z_{k}^{q_{k}}+g_{k}(z_0,\dots,z_{k-1}),
  \end{align*}
where the $g_j(z_0,\dots,z_{j-1})$ are as before.

We again successively substitute each of these equations into the next
to eliminate the variables $z_2,z_3, \dots,z_k$.  To be precise, we
first make these substitutions to put the equations in the form:
\begin{align*}
z_2&=a_{1}z_1^{q_1}+a_{0}z_0^{p_1}\\
z_3&=f_2(z_0,z_1)\\
\dots&\quad\dots\quad \dots\\
z_k&=f_{k-1}(z_0,z_1)\\
z^n&=f_{k}(z_0,z_{1}).
  \end{align*}
  Recall our notation $f=f_{k}$. In terms of new coordinates, $x=z_0$,
  $y=z_1$, $z$, $ Z_2:=z_2-a_{1}z_1^{q_1}+a_{0}z_0^{p_1}$, $ \dots$, $
  Z_{k}:=z_k-f_{k-1}(z_0,z_1)$, these equations become
$$
Z_2=Z_3=\dots=Z_{k}=0; \quad z^n=f(x,y).
$$

We are now ready to prove the main result of this section.
\begin{theorem}
  Let $(X,o)$ be a hypersurface singularity at the origin given by an
  equation in the form $z^n+g(x,y)=0$ with homology sphere link. Then
  the Milnor Fiber Conjecture is true for $(X,o)$.
\end{theorem}
\begin{proof}
  Suppose that we have a splice decomposition corresponding to the
  following decomposition of our splice diagram as the splice of two
  diagrams:
$$\splicediag{16}{24}{
  \Circ\lineto[r]^(.75){p_1}&\Circ\lineto[r]^(.25){n}^(.75){p_2}
  \lineto[d]^(.25){q_1}&\Circ\lineto[d]^(.25){q_2}\lineto[r]^(.25)n&
  \dotto[r]&\lineto[r]^(.75){p_r}&\Circ\lineto[d]^(.25){q_r}\ar[r]^(.25)n&
  \quad&\Circ\ar[l]_(.25){p_{r+1}}\lineto[d]^(.25){q_{r+1}}\lineto[r]^(.25)n&
  \dotto[r]&\lineto[r]^(.75){p_k}&\Circ\lineto[r]^(.25)n\lineto[d]^(.25){q_k}&
  \Circ\\
  &\Circ&\Circ&&&\Circ&&\Circ&&&\Circ}
$$
We wish to show that the Milnor fiber $F$ for $z^n=g(x,y)$ is
obtained by the construction $F_1^o\cup_{N_1}(G_1\times
G_2)\cup_{N_2}F_2^o$ of Conjecture \ref{conj:milnor fiber}, where $F_1$
and $F_2$ are Milnor fibers for the two splice components, $G_1$
and $G_2$ are fibers in the links of the two splice components for the
knots along which we splice, and $F_i^o$ is the result of removing a
tubular neighborhood of a properly embedded $G_i$ in $F_i$.

In \cite{neumann-bams} (see also \cite{kauffman-neumann}) it is shown
that the Milnor fiber $F$ is obtained by taking a Milnor fiber
$G\subset S^3=\partial D^4$ for $g$, pushing it inside $D^4$ so that
it is properly embedded (that is, $\partial G= G\cap\partial D^4$),
and then taking the $n$--fold branched cyclic cover of $D^4$, branched
along this embedding of $G$.

We need to understand the placement of $G$ with respect to the
decomposition of $D^4$ of equation (\ref{eq:3}). On taking the n-fold
branched cover we will see that we get the desired decomposition of
$F$.

According to \cite{eisenbud-neumann} the fiber $G$ decomposes
according to the splice diagram into $q_{r+1}\dots q_k$ parallel
copies of the Milnor fiber $G_1'$ of the plane curve given by
$$\splicediag{16}{24}{
  \Circ\lineto[r]^(.75){p_1}&\Circ\lineto[r]^(.25){1}^(.75){p_2}
  \lineto[d]^(.25){q_1}&\Circ\lineto[d]^(.25){q_2}\lineto[r]^(.25)1&
  \dotto[r]&\lineto[r]^(.75){p_r}&\Circ\lineto[d]^(.25){q_r}\ar[r]^(.25)1&
  \\
  &\Circ&\Circ&&&\Circ& } $$
and one copy of the Milnor fiber of the
plane curve corresponding to
$$\splicediag{16}{24}{
  \quad\Circ&\Circ\lineto[l]_(.25){p_{r+1}}
  \lineto[d]^(.25){q_{r+1}}\lineto[r]^(.25)1&
  \dotto[r]&\lineto[r]^(.75){p_k}&\Circ\ar[r]^(.25)1\lineto[d]^(.25){q_k}&
  \\
  &\Circ&&&\Circ}
$$
punctured $q_{r+1}\dots q_k$ times.

We can position $G$ with respect to the decomposition of equation
(\ref{eq:3}) so that it lies completely in $(G_1'\times D^2)\cup
(F_2')^o$. It then intersects $(G_1'\times D^2)$ in $q_{r+1}\dots q_k$
parallel copies of $G_1'$. Its intersection with $(F_2')^o$ is
obtained as follows.  First make the
fiber $G'$ of the knot represented by the right arrowhead of the splice
diagram:
$$\splicediag{16}{24}{
    &\Circ\ar[l]_(.25){p_{r+1}}\lineto[r]^(.25){1}^(.75){p_{r+2}}
    \lineto[d]^(.25){q_{r+1}}&\Circ\lineto[d]^(.25){q_{r+2}}\lineto[r]^(.25)1&
    \dotto[r]&\lineto[r]^(.75){p_k}&\Circ\lineto[d]^(.25){q_k}\ar
[r]^(.25)1&\\
    &\Circ&\Circ&&&\Circ}
  $$
  properly embedded in $D^4$ and transverse to the properly
  embedded version of the fiber $D^2$ of the unknot represented by the
  left arrowhead. Then remove the tubular neighborhood of the latter.
  Using \cite{neumann-bams}, the $n$--fold cyclic cover of $D^4$ along
  $G'$ is the Milnor fiber for the surface singularity with diagram:
$$\splicediag{16}{24}{
   \Circ &\Circ\lineto[l]_(.25){p_{r+1}}\lineto[r]^(.25){n}^(.75){p_{r+2}}
    \lineto[d]^(.25){q_{r+1}}&\Circ\lineto[d]^(.25){q_{r+2}}\lineto[r]^(.25)n&
    \dotto[r]&\lineto[r]^(.75){p_k}&\Circ\lineto[d]^(.25){q_k}\lineto
[r]^(.25)n&\Circ\\
    &\Circ&\Circ&&&\Circ}
  $$
Moreover, the embedded $D^2\subset D^4$ lifts in this cover to copy of
the fiber for the knot represented by the left-most vertex.

It follows that the decomposition of equation (\ref{eq:3}) lifts to
give the desired decomposition of $F$, as desired.
\end{proof}

In the context of the above result it is worth mentioning that
N\'emethi and Mendris recently showed \cite{nemethi-mendris} that for
a singularity $z^n=f(x,y)$ with homology sphere link (even rational
homology sphere link) the Milnor fibration is topologically determined
by the link of the singularity.
\begin{remark}\label{collin-saveliev}
  The results of this section give a proof of the Casson Invariant
  Conjecture (CIC) for these examples, also proven in
  \cite{neumann-wahl90, collin-saveliev,
  nemethi-nicolaescu3}. Saveliev and Collin
  \cite{collin-saveliev}, using
  equivariant Casson invariant, give an iterative generalization of
  these examples but their approach implies more:
  Let $\Delta$ be any splice diagram satisfying the semigroup
  condition and $w$ a leaf of $\Delta$. We allow, as in this section,
  the weight on the edge to $w$ to be $1$.  For $n\in\N$ let
  $\Delta_n(w)$ be the diagram obtained by multiplying the weight
  furthest from $w$ on each edge by $n$. We assume $n$ is chosen
  coprime to all the unchanged weights at each node, so $\Delta_n(w)$
  is again a splice diagram. Then if CIC
  is valid for splice type singularities for $\Delta$, then the same
  holds for $\Delta_n(w)$.
\end{remark}

\section{Appendix: Splicing and plumbing}\label{sec:splicing}

In this appendix we recall the classification of $\Z$--homology
sphere singularity links in terms of splice diagrams and describe how
to recover a resolution diagram from the splice diagram.

We start with Seifert fibered manifolds. For the following results see
\cite{neumann-raymond}. Let $\Sigma$ be a Seifert fibered homology
3--sphere other than $S^3$.  Then it has at least $3$ singular fibers
and the degrees $p_1, \dots ,p_r$ of these singular fibers are pairwise
coprime. Conversely, given a set $\{p_1, \dots, p_r\}$ of pairwise
coprime integers $p_i>1$ with $r\ge3$, there is a unique Seifert
fibered homology sphere $\Sigma(p_1,\dots,p_r)$ up to orientation with
these singular fiber degrees. Moreover, $\Sigma(p_1,\dots,p_r)$ has a
unique orientation for which it is a singularity link, so we give it
this orientation.  It is, in fact, the link of the \BCI
%\Brieskorn{} complete intersection singularity
$$V(p_1,\dots,p_r):=
\{(z_1,\dots,z_r)\in\C^n:a_{i1}z_1^{p_1}+\dots+a_{ir}z_r^{p_r}=0
\text{ for } i=1,\dots, r-2\},$$
for a sufficiently general matrix
$(a_{ij})$ of coefficients. By Hamm \cite{hamm}, ``sufficiently
general'' means that all $(r-2)\times(r-2)$ minors should be
non-singular.

We represent the homology sphere $\Sigma(p_1,\dots,p_r)$ by the
\emph{splice diagram\/}:
$$
\splicediag{10}{36}{
\Circ&{\raise5pt\hbox{$\bf \dots\dots$}}&\Circ\\
&\Circ\lineto[dl]^(.25){p_r}\lineto[ul]_(.25){p_{r-1}}
\lineto[ur]^(.25){p_2}\lineto[dr]_(.25){p_1}\\
\Circ&&\Circ}$$
Each of the singular fibers of $\Sigma(p_1,\dots,p_r)$ represents a
knot in $\Sigma(p_1,\dots,p_r)$ which we represent in a splice diagram
by adding an arrowhead to the corresponding edge. Thus
$$\splicediag{12}{24}{
&&\\
\Circ&\Circ\ar[dr]_(.25){2}\ar[ur]^(.25){3}
\lineto[l]_(.25){5}\\
&&
}$$
represents the link in $\Sigma(2,3,5)$ consisting of the
degree $2$ and $3$ singular fibers. Non-singular fibers are
represented by adding new arrows at the central vertex weighted by 1,
so
$$\splicediag{12}{24}{
\Circ&&\Circ\\
&\Circ\lineto[dr]_(.25){2}\lineto[ur]^(.25){3}
\lineto[ul]_(.25){5}\ar[dl]^(.25){1}\\
&&\Circ
}$$
 represents the knot in $\Sigma(2,3,5)$ consisting of one
non-singular fiber.

There are Seifert fibrations of the $3$--sphere with $2$ or less
singular fibers. For instance, $S^3$ can be fibered by copies of the
$(p,q)$ torus knot, with one $p$--fold singular fiber and one $q$--fold
singular fiber, so the splice diagram
$$\splicediag{12}{24}{
&&\Circ\\
&\Circ\lineto[dr]_(.25){p}\lineto[ur]^(.25){q}
\ar[l]_(.25){1}\\
&&\Circ
}$$
is the diagram for the $(p,q)$ torus knot in $S^3$. Similarly
$$\splicediag{12}{24}{
&&\\
\Circ&\Circ\ar[dr]_(.25){1}\ar[ur]^(.25){1}
\lineto[l]_(.25){q}\\
&&
}$$
represents a pair of parallel $(1,q)$ torus knots
(unknotted curves which link each other $q$ times).

If $K_1\subset \Sigma_1$ is a knot in a homology sphere and
$K_2\subset \Sigma_2$ is another, then we form the \emph{splice} of
$\Sigma_1$ to $\Sigma_2$ along $K_1$ and $K_2$ as follows.  Let $N_i$
be a closed tubular neighborhood of $K_i$ in $\Sigma_i$ for $i=1,2$
and let $\Sigma_i'$ be the result of removing its interior, so
$\partial \Sigma_i'= T^2$. The splice is the manifold
$$\Sigma=\Sigma_1'\cup_{T^2}\Sigma_2'\,,$$
where the gluing matches
meridian in $\Sigma_1$ to longitude in $\Sigma_2$ and vice versa.
(``Meridian'' and ``longitude'' in $\Sigma_1'$ are the simple curves
in $\partial \Sigma_1'=T^2$ that are null-homologous respectively in
the removed solid torus $N_1$ or in $\Sigma_1'$.) We denote the splice
by
$$\Sigma=\Sigma_1~\raise5pt\hbox{$\underline{K_1\quad K_2}$}~\Sigma_2\,.$$

We represent splicing in terms of splice diagrams by gluing the
diagrams at the arrowheads that represent the knots along which we are
splicing.  For instance,
$$\splicediag{12}{30}{
\Circ&&&\Circ\\
&\Circ\lineto[ul]_(.25){2}\lineto[dl]^(.25)3
&\Circ\lineto[dr]_(.25){2}\lineto[ur]^(.25){3}
\lineto[l]_(.2){7}_(.8){7}\\
\Circ&&&\Circ
}$$
represents the splice of two copies of $\Sigma(2,3,7)$
along the knots represented by the degree $7$ fibers.

By \cite{eisenbud-neumann}, the splice diagrams that classify homology
sphere singularity links are precisely the splice diagrams with
pairwise coprime positive weights around each node and with positive
edge determinants (recall that the \emph{edge determinant} is the
product of the two weights on the edge minus the product of the
weights adjacent to the edge).

The splice diagram can be computed very easily from a resolution
diagram for the singularity.  We describe this in detail in the
appendix of \cite{neumann-wahl10} so we will not repeat it here.
Briefly, the splice diagram is obtained from the dual resolution graph
for the singularity by replacing each string in the resolution graph
by a single edge (ie, we eliminate vertices of valence 2);
each splice diagram weight is the absolute value of the determinant of
the intersection matrix for the subgraph of the resolution graph cut
off at the corresponding node in the direction of the corresponding
edge.

We will also need an ``unreduced'' version of the splice diagram: The
\emph{maximal splice diagram} is the version of the splice diagram we
get from the resolution graph if we do not first eliminate vertices of
valency $2$, and we include edge weights at \emph{all} vertices ---
also the leaves.  For example, the
resolution graph
$$
\xymatrix@R=6pt@C=24pt@M=0pt@W=0pt@H=0pt{
\\
\overtag{\Circ}{-2}{8pt}&&&&\overtag{\Circ}{-2}{8pt}\\
&\overtag{\Circ}{-1}{8pt}\lineto[ul]\lineto[dl]\lineto[r]&
\overtag{\Circ}{-17}{8pt}&\overtag{\Circ}{-1}{8pt}\lineto[ur]\lineto[dr]\lineto[l]&\\
\overtag{\Circ}{-3}{8pt}&&&&\overtag{\Circ}{-3}{8pt}\lineto[r]&\overtag{\Circ}{-2}{8pt}}
$$
gives maximal splice
diagram and splice diagram
$$
\xymatrix@R=6pt@C=24pt@M=0pt@W=0pt@H=0pt{
    \Circ&&&&\Circ\\
  &\Circ\lineto[ul]_(.25){2}_(.75){11}
\lineto[dl]^(.25){3}^(.75)5
\lineto[r]^(.25)7^(.75)1&
\Circ&\Circ\lineto[ur]^(.25)2^(.75){28}\lineto[dr]_(.25)5_(.75)9
\lineto[l]_(.25){11}_(.75)1&\\
\Circ&&&&\Circ\lineto[r]_(.25)2_(.75)5&\Circ}
%$$
%and the splice diagram
%$$
\splicediag{6}{30}{
&\Circ&&&\Circ\\
\qquad\text{and}&&\Circ\lineto[ul]_(.25){2}\lineto[dl]^(.25)3
&\Circ\lineto[dr]_(.25){5}\lineto[ur]^(.25)2
\lineto[l]_(.2){11}_(.8){7}\\
&\Circ&&&\Circ
}\qquad$$
respectively  (this is Example 12.1 in \cite{neumann-wahl10}).

An algorithm to recover the resolution diagram from the splice diagram
is given in \cite{eisenbud-neumann}. Here we describe an easier method
that arose from conversations with Paul Norbury (developed
independently by Pierrette Cassou-Nogues \cite{cassou-nogues}, whose
terminology of ``maximal splice diagram'' we have adopted---we called
it ``adjoint diagram'').  

To compute the resolution graph from the splice diagram we will give
algorithms to:
\begin{itemize}
\item compute the maximal splice diagram from the splice diagram, and
\item compute the resolution graph from the maximal splice diagram.
\end{itemize}

We will need the following properties of the maximal splice diagram,
which are proved in greater generality in section 12 of
\cite{neumann-wahl10} (Theorem 12.2 and Lemma 12.5).
\begin{theorem}\label{th:props}
  {\rm(1)}\qua For any pair of vertices $v$ and $w$ of the maximal diagram
  let $\ell_{vw}$ be the product of the weights adjacent to, but not on,
  the shortest path from $v$ to $w$ in $\Delta'$. Then the matrix
  $L:=(\ell_{vw})$ is the inverse matrix of $-A(\resgraph)$.

  {\rm(2)}\qua   Every edge determinant for the maximal splice diagram is
  $1$.

  {\rm(3)}\qua  The edge-weight adjacent to a leaf $v$ of the maximal splice
  diagram is equal to $\lceil a/b\rceil$ where $a$ is the product of
  edge-weights adjacent to and just beyond the nearest node to $v$ and
  $b$ is the remaining weight adjacent to that node.\qed
 \end{theorem}
 We remark that part (3) is valid also for the valency $2$ vertices
 between the leaf and its nearest node. For example, for the
 right-most leaf of the above
 example $5=\lceil 22/5\rceil=\lceil9/2\rceil$.

\subsection{Maximal splice diagram from splice diagram}\label{subsec:9.1}
We describe how to recover the string of vertices and weights of the
maximal splice diagram between any two vertices of a splice diagram.
Suppose first both vertices are nodes with weights as follows,
$$\splicediag{8}{40}{ &&&\\&&&\\
  \Vdots&\Circ\lineto[uul]_(.35){a_1}
  \lineto[ul]^(.45){a_2}\lineto[dl]^(.35){a_r}
  \lineto[r]^(.25)b^(.75)c&\Circ
  \lineto[uur]^(.35){d_1}\lineto[ur]_(.45){d_2}\lineto[dr]_(.35){d_s}
  &\Vdots\\
  &&&\\~ }
$$
and put $a=\prod_1^r a_i$, $d=\prod_1^s d_j$. If one of the
vertices (say the right one) is a leaf instead of a node then we put
$d=1$.  The desired string of vertices and weights between our two
nodes will only depend on $a,b,c,d$, so we replace the above diagram
by:
$$\splicediag{6}{18}{
\lineto[r]^(.7)a&\Circ\lineto[rr]^(.2)b^(.8)c&&\Circ\lineto[r]^(.3)d&
}$$

Consider the following infinite linear graph:
$$\splicediag{6}{18}{
&\dotto[l]\lineto[r]^(.7)1&
\Circ\lineto[rr]^(.15)3^(.85)1&&
\Circ\lineto[rr]^(.15)2^(.85)1&&
%\Circ\lineto[rr]^(.25)1^(.75)1&&
\Circ\lineto[rr]^(.15)1^(.85)2&& \Circ\lineto[rr]^(.15)1^(.85)3&&
\Circ&\lineto[l]_(.7)1\dotto[r]& }$$
We are going to refine this by
adding vertices on this line until our vertices
$\splicediag{6}{15}{\lineto[r]^(.7)a&\Circ\lineto[r]^(.3)b&}$ and
$\splicediag{6}{15}{\lineto[r]^(.7)c&\Circ\lineto[r]^(.3)d&}$ appear
on it. Vertices
$\splicediag{6}{15}{\lineto[r]^(.7)x&\Circ\lineto[r]^(.3)y&}$ are
ordered along the line by size of $x/y$. Thus such a vertex either is
already a vertex of the linear graph, or it falls on an existing edge.
In the latter case we subdivide the edge as follows:
$$\splicediag{6}{24}{
  \lineto[r]^(.7)\alpha&\Circ\lineto[rr]^(.15)\beta^(.85)\gamma&&
  \Circ\lineto[r]^(.3)\delta& \\\\} \quad\mapsto\quad
\splicediag{6}{20}{
  \lineto[r]^(.7)\alpha&\Circ\lineto[rr]^(.15)\beta^(.75){\alpha+\gamma}&&
  \Circ\lineto[rr]^(.25){\beta+\delta}^(.85)\gamma&&
  \Circ\lineto[r]^(.3)\delta&\\\\}$$
We repeat this process until
both our desired vertices appear, and then the portion of the linear
graph between them is what we were seeking.

For example, suppose our initial splice diagram is: 
$$\splicediag{6}{24}{
\Circ\lineto[dr]^(.7)2&&&&\Circ\lineto[dl]_(.7)2\\
&\Circ\lineto[rr]^(.2){7}^(.8){11}&&
\Circ\\
\Circ\lineto[ur]_(.7)5&&&&\Circ\lineto[ul]^(.7)3\\
}$$
To create the string for the middle edge we start with: 
$$\splicediag{6}{18}{
  \lineto[r]^(.6){10}&\Circ\lineto[rr]^(.2){7}^(.8){11}&&
  \Circ\lineto[r]^(.4){6}& }$$
and apply the above procedure. We mark the
positions of these vertices, until they are found, by
\,\,{\scriptsize$\vee$} .
\begin{gather*}
\splicediag{6}{80}{\\
\lineto[r]^(.9)1&\Circ\lineto[rr]^(.05){1}
^(.35){{}^{10}\vee^7}^(.65){{}^{11}\vee^6}
^(.95)2&&\Circ\lineto[r]^(.1){1}&
}\\
\splicediag{6}{52}{
\lineto[r]^(.88)1&\Circ\lineto[rr]^(.06){1}
^(.6){{}^{10}\vee^7}^(.94)3&&\Circ\lineto[rr]^(.06)2
^(.6){{}^{11}\vee^6}
^(.94)2&&\Circ\lineto[r]^(.12){1}&
}\\
\splicediag{6}{25.5}{
\lineto[r]^(.7)1&\Circ\lineto[rr]^(.15){1}^(.85)4
&&\Circ\lineto[rrr]^(.1)3
^(.6){{}^{10}\vee^7}^(.9)3&&&\Circ\lineto[rr]^(.15)2^(.85)5&&
\Circ\lineto[rrr]^(.1)3
^(.6){{}^{11}\vee^6}
^(.9)2&&&\Circ\lineto[r]^(.3){1}&
}\\
\splicediag{6}{18.7}{
\lineto[r]^(.7)1&\Circ\lineto[rr]^(.15){1}^(.85)4
&&\Circ\lineto[rr]^(.15)3^(.85)7
&&\Circ\lineto[rrr]^(.1)5
^(.5){{}^{10}\vee^7}^(.9)3&&&\Circ\lineto[rr]^(.15)2^(.85)5&&
\Circ\lineto[rr]^(.15)3^(.85)7
&&\Circ\lineto[rrr]^(.1)4
^(.6){{}^{11}\vee^6}
^(.9)2&&&\Circ\lineto[r]^(.3){1}&
}\\
\splicediag{6}{15.7}{
\lineto[r]^(.6)1&\Circ\lineto[rr]^(.2){1}^(.8)4
&&\Circ\lineto[rr]^(.2)3^(.8)7
&&\Circ\lineto[rr]^(.2)5^(.8){10}
&&\Dot\lineto[rr]
^(.2)7^(.8)3&&\Circ\lineto[rr]^(.2)2^(.8)5&&
\Circ\lineto[rr]^(.2)3^(.8)7
&&\Circ\lineto[rr]^(.2)4^(.8)9
&&\Circ\lineto[rrr]^(.15)5
^(.5){{}^{11}\vee^6}
^(.85)2&&&\Circ\lineto[r]^(.4){1}&
}\\
\splicediag{6}{15}{
\lineto[r]^(.6)1&\Circ\lineto[rr]^(.2){1}^(.8)4
&&\Circ\lineto[rr]^(.2)3^(.8)7
&&\Circ\lineto[rr]^(.2)5^(.8){10}
&&\Dot\lineto[rr]
^(.2)7^(.8)3&&\Circ\lineto[rr]^(.2)2^(.8)5&&
\Circ\lineto[rr]^(.2)3^(.8)7
&&\Circ\lineto[rr]^(.2)4^(.8)9
&&\Circ\lineto[rr]^(.2)5^(.8){11}
&&\Dot\lineto[rr]^(.2)6
^(.8)2&&\Circ\lineto[r]^(.4){1}&
\\ \\}
\end{gather*}
Thus the final string is:
$$\splicediag{6}{25}{\\
\lineto[r]^(.75){10}
&\Dot\lineto[rr]^(.1)7^(.9)3
&&\Circ\lineto[rr]^(.1)2^(.9)5
&&\Circ\lineto[rr]^(.1)3^(.9)7
&&\Circ\lineto[rr]^(.1)4^(.9)9
&&\Circ\lineto[rr]^(.1)5^(.85){11}
&&\Dot\lineto[r]^(.2)6&\\\\
}$$
Similarly, the $3$--weighted edge expands from 
$~\splicediag{6}{12}{\lineto[r]^(.7){22}&
\Circ\lineto[rr]^(.2)3^(.8)8&&\Circ\lineto[r]^(.3)1&}~$
as follows:
\begin{gather*}
\splicediag{6}{80}{\\
\lineto[r]^(.9)7&\Circ\lineto[rr]^(.05){1}
^(.4){{}^{22}\vee^3}^(.95)8&&\Dot\lineto[r]^(.1){1}&
}\\
\splicediag{6}{52}{
\lineto[r]^(.88)7&\Circ\lineto[rr]^(.05){1}^(.5){{}^{22}\vee^3}^(.95){15}&&
\Circ\lineto[r]^(.1)2^(.9)8&\Dot\lineto[r]^(.1)1&
}\\
\splicediag{6}{52}{
\lineto[r]^(.88)7&\Circ\lineto[r]^(.1){1}^(.9){22}&
\Dot\lineto[r]^(.1)3^(.9){15}&
\Circ\lineto[r]^(.1)2^(.9)8&\Dot\lineto[r]^(.1)1&
}
\end{gather*}
A shortcut is available in the above procedure: to compute the string
for the central edge we did not need to create the $(4,3)$-- and
$(7,5)$--vertices, since the edge determinant on the edge from
$(10,7)$ to $(3,2)$ is already 1.

With this comment, the other three edges are immediate and the maximal
splice diagram is:
$$\splicediag{6}{16}{
\Circ\lineto[drr]^(.2){18}^(.8)2&&&&&&&&&&&&
&&\Circ\lineto[dll]_(.2){17}_(.8)2\\
&&\Circ\lineto[rr]^(.2)7^(.8)3
&&\Circ\lineto[rr]^(.2)2^(.8)5
&&\Circ\lineto[rr]^(.2)3^(.8)7
&&\Circ\lineto[rr]^(.2)4^(.8)9
&&\Circ\lineto[rr]^(.2)5^(.8){11}
&&\Circ\\
\Circ\lineto[urr]_(.2)3_(.8)5&&&&&&&&&&&&
&&\Circ\lineto[ull]^(.8)3^(.2){15}
\lineto[rr]^(.2)2^(.8)8&&\Circ\\
}$$

\subsection{Resolution graph from maximal splice diagram} \label{subsec:9.2}
We must recover the self-intersection weights
$e_v:=a_{vv}$ at vertices.  The matrix equation $LA(\resgraph)=-I$ gives
equations that will do this. We use the notation $w\hbox{--}v$ to mean
vertices $w$ and $v$ are connected by an edge. Then
%$$e_v=\frac{-1}{\ell_{vv}}
%\left(1+\sum_{\{w:w\lower1pt\hbox{--}v\}}{\ell_{vw}}\right)\,.$$
%Alternatively, if
for any vertex $w'$ adjacent to $v$, the $vw'$ entry of this matrix
equation gives:
$$e_v=\frac{-1}{\ell_{vw'}}\left(
\sum_{{\{w:w\lower1pt\hbox{--}v\}}
}\ell_{ww'}
\right)
$$
Note that the product of the weights just beyond $w'$ from $v$ cancel
in this formula, so they may be replaced by $1$ for the calculation.
For example, for
the above maximal splice diagram we get the resolution graph:
$$\splicediag{6}{16}{\\
\overtag\Circ{-2}{8pt}&&&&&&&&&&&&&&\overtag\Circ{-2}{8pt}\\
%\lineto[r]
&&\overtag\Circ{-1}{8pt}\lineto[rr]\lineto[ull]\lineto[dll]
&&\overtag\Circ{-5}{8pt}\lineto[rr]
&&\overtag\Circ{-2}{8pt}\lineto[rr]
&&\overtag\Circ{-2}{8pt}\lineto[rr]
&&\overtag\Circ{-2}{8pt}\lineto[rr]
&&\overtag\Circ{-2}{8pt}\lineto[urr]\lineto[drr]\\
\overtag\Circ{-5}{8pt}&&&&&&&&&&&&&&
\overtag\Circ{-2}{8pt}\lineto[rr]&&\overtag\Circ{-2}{8pt}\\\\
}$$

\subsection{Proof of the procedure of \ref{subsec:9.1}}\label{subsec:proof} 
The procedure in \cite{eisenbud-neumann} implies that the string of
the maximal splice diagram between the two
vertices in question only depends on $a,b,c,d$. Consider the
resolution graph
\begin{equation}
  \label{eq:1}
  \splicediag{6}{10}{
\overtag\Circ{-1}{6pt}\lineto[rr]&&
\overtag\Circ{-2}{6pt}\lineto[r]&\dotto[r]&
&\dotto[r]&
\overtag\Circ{-2}{6pt}\lineto[rr]&&
\overtag\Circ{-3}{6pt}\lineto[rr]&&
\overtag\Circ{-2}{6pt}\lineto[rr]&&
\overtag\Circ{-2}{6pt}\lineto[rr]&&
\overtag\Circ{-2}{6pt}\dotto[r]&
&\dotto[r]&
\overtag\Circ{-2}{6pt}\lineto[rr]&&
\overtag\Circ{-1}{6pt}
}
\end{equation}
with associated maximal splice
diagram:
\begin{equation}
  \label{eq:2}
  \splicediag{6}{10}{
\Circ\lineto[rr]^(.25)s^(.75)1&&
\Circ\lineto[r]^(.6){s-1}&\dotto[r]&
&\dotto[r]&
\Circ\lineto[rr]^(.25)2^(.75)1&&
\Circ\lineto[rr]^(.25)1^(.75)2&&
\Circ\lineto[rr]^(.25)1^(.75)3&&
\Circ\lineto[rr]^(.25)1^(.75)4&&
\Circ\dotto[r]&
&\dotto[r]&
\Circ\lineto[rr]^(.25)1^(.75)t&&
\Circ
}
\end{equation}
This is a piece of the infinite linear graph we used above, and we
choose $s$ and $t$ large enough that our desired vertices will lie
in this piece. Now we repeatedly blow up on edges of the linear
resolution graph. An easy calculation shows that blowing up on an
edge:
$$\splicediag{6}{20}{
\lineto[r]&\overtag\Circ{e_1}{8pt}\lineto[rr]&&
\overtag\Circ{e_2}{8pt}\lineto[r]
&&\mapsto&
\lineto[r]&\overtag\Circ{e_1-1}{8pt}\lineto[rr]&&
\overtag\Circ{-1}{8pt}\lineto[rr]&&
\overtag\Circ{e_2-1}{8pt}\lineto[r]&}$$
has the effect:
$$\splicediag{6}{20}{
  \lineto[r]^(.6)\alpha&\Circ\lineto[rr]^(.2)\beta^(.8)\gamma&&
  \Circ\lineto[r]^(.4)\delta& &\mapsto&
  \lineto[r]^(.6)\alpha&\Circ\lineto[rr]^(.2)\beta^(.75){\alpha+\gamma}&&
  \Circ\lineto[rr]^(.25){\beta+\delta}^(.8)\gamma&&
  \Circ\lineto[r]^(.4)\delta&\\\\}$$
on the associated maximal splice diagram. Thus we need only show that
our desired vertices eventually appear in this procedure. But this is
a standard fact about Farey sequences (alternatively, one can observe
that we are describing the standard procedure to resolve the
plane curve singularity $(x^a+y^b)(x^c+y^d)$).

This same argument applies to see how to fill in the maximal splice
diagram between a node and a leaf, even if the edge weight at the leaf
is unknown. The leaf will be the rightmost vertex of the above string
(\ref{eq:2}) with $t$ chosen as small as possible to accommodate our
desired vertex
$\splicediag{6}{15}{\lineto[r]^(.5)a&\Circ\lineto[r]^(.5)b&}$. Thus,
the $t$ that we choose is $\lceil a/b\rceil$ (if $t=1$ the initial
resolution string (\ref{eq:1}) is
\,\,$\phantom{\Bigl(}\splicediag{6}{10}{
  \overtag\Circ{-1}{6pt}\lineto[rr]&&
  \overtag\Circ{-2}{6pt}\lineto[r]&\dotto[r]& &\dotto[r]&
%\overtag\Circ{-2}{6pt}\lineto[rr]&&
\overtag\Circ{-2}{6pt}}$\,\,).

Note that this blow-up procedure gives an alternative way to compute
the self-intersection weights along the string, making the calculation
of subsection \ref{subsec:9.2} only necessary at nodes.

\end{document}